\documentclass[article,11pt]{amsart}
\usepackage{amsmath,mathrsfs,bm,amssymb,color}%suzuki
\usepackage{graphicx}

\usepackage{caption}
\usepackage{subcaption}

\newtheorem{theorem}{Theorem}[section]
\newtheorem{proposition}[theorem]{Proposition}
\newtheorem{lemma}[theorem]{Lemma}
\newtheorem{corollary}[theorem]{Corollary}
\newtheorem{definition}[theorem]{Definition}

\newtheorem{remark}[theorem]{Remark}

\newtheorem{assumption}[theorem]{Assumption}
\newcommand{\bl}[1]{\begin{lemma}\label{#1}}
\newcommand{\el}{\end{lemma}}
\newcommand{\bas}[1]{\begin{assumption}\label{#1}}
\newcommand{\eas}{\end{assumption}}
\newcommand{\bc}[1]{\begin{corollary}\label{#1}}
\newcommand{\ec}{\end{corollary}}
\newcommand{\bt}[1]{\begin{theorem}\label{#1}}
\newcommand{\et}{\end{theorem}}
\newcommand{\bp}[1]{\begin{proposition}\label{#1}}
\newcommand{\ep}{\end{proposition}}
\newcommand{\br}[1]{\begin{remark}\label{#1}}
\newcommand{\er}{\end{remark}}
\newcommand{\bdf}[1]{\begin{definition}\label{\rm #1}}
\newcommand{\edf}{\end{definition}}
\newcommand{\bi}{\begin{description}}
\newcommand{\ei}{\end{description} }

\newcommand{\eq}[1]{\begin{equation}\label{#1}}
\newcommand{\en}{\end{equation}}
\newcommand{\eqst}{\begin{equation*}}
\newcommand{\enst}{\end{equation*}}
\newcommand{\eqn}{\begin{eqnarray*}}
\newcommand{\enn}{\end{eqnarray*}}
\newcommand{\eqnn}{\begin{eqnarray}}
\newcommand{\ennn}{\end{eqnarray}}
\newcommand{\ex}{\mathbb{E}}

\DeclareMathOperator{\Ai}{{\rm Ai}}
\DeclareMathOperator{\Aip}{{\rm Ai^{\prime}}}

\DeclareMathOperator{\Aik}{{{\rm Ai}_k}}

\DeclareMathOperator{\rA}{{\rm A}}

\newcommand{\ap}{\alpha}
\newcommand{\BT}{\beta}

\newcommand{\EP}{\epsilon}
\newcommand{\GA}{\Gamma}
\newcommand{\ga}{\gamma}

\newcommand{\del}{\delta}

\newcommand{\s}{\sigma}
\newcommand{\LA}{\Lambda}
\newcommand{\la}{\lambda}

\newcommand{\ud}{\, \mathrm{d}}

\newcommand{\NN}{\mathbb{N}}
\newcommand{\ZZ}{\mathbb{Z}}
\newcommand{\RR}{\mathbb{R}}

\newcommand{\cF}{\mathcal{F}}
\newcommand{\cG}{\mathcal{G}}

\newcommand{\cL}{\mathcal{L}}

\newcommand{\cP}{\mathcal{P}}
\newcommand{\cQ}{\mathcal{Q}}

\newcommand{\rF}{{\rm F}}
\newcommand{\sF}{\mathscr{F}}

\newcommand{\lkkk }{\left[}
\newcommand{\rkkk}{\right]}
\newcommand{\lk}{\left(}
\newcommand{\rk}{\right)}
\newcommand{\lkk}{\left\{}
\newcommand{\rkk}{\right\}}
\newcommand{\lka }{\left\langle}
\newcommand{\rka}{\right\rangle}

\newcommand{\non}{\nonumber}

\providecommand{\Tr}{\mathop{\mathrm{Tr}}\nolimits}
\providecommand{\Spec}{\mathop{\mathrm{Spec}}\nolimits}
\providecommand{\pro}[1]{(#1_t)_{t\geq0}}

\makeatletter
\@addtoreset{equation}{section}
\makeatother

\title{Spectral Properties of the Massless Relativistic Quartic Oscillator}
\author{Samuel O. Durugo and J\'ozsef L\H{o}rinczi*}
\address{\it \small Department of Mathematical Sciences, Loughborough University, Loughborough LE11 3TU, United Kingdom}
\email{\small {\tt  S.O.Durugo4@lboro.ac.uk, J.Lorinczi@lboro.ac.uk}}

\date{}

\begin{document}
\thanks{*corresponding author}
\maketitle
\setlength{\baselineskip}{14pt}

\begin{abstract}
\noindent
An explicit solution of the spectral problem of the non-local Schr\"odinger operator obtained as the sum of the
square root of the Laplacian and a quartic potential in one dimension is presented. The eigenvalues are obtained
as zeroes of special functions related to the fourth order Airy function, and closed formulae for the Fourier
transform of the eigenfunctions are derived. These representations allow to derive further spectral properties
such as estimates of spectral gaps, heat trace and the asymptotic distribution of eigenvalues, as well as a
detailed analysis of the eigenfunctions. A subtle spectral effect is observed which manifests in an
exponentially tight approximation of the spectrum by the zeroes of the dominating term in the Fourier
representation of the eigenfunctions and its derivative.


\medskip
\noindent
2010 {\it MS Classification}: Primary 35P15, 47G30, 81Q05; Secondary 35P05,  60G52, 81Q10 \\
\end{abstract}
\bigskip

%\newpage
\section{Introduction}
In this paper our aim is to obtain explicit formulae on the spectrum and eigenfunctions of the specific
non-local Schr\"odinger operator
\begin{equation}
\label{quart}
H = \sqrt{-\frac{\ud^{2}}{\ud x^{2}}} + x^4,
\end{equation}
defined as the sum in an appropriate sense of the one-dimensional fractional Laplacian of index $\frac{1}{2}$
and the potential $V(x) = x^4$ acting as a multiplication operator. While for the classical Schr\"odinger
operator featuring the Laplacian and the same potential there are no closed-form solutions of the eigenvalue
problem, this can be obtained for the relativistic Schr\"odinger operator above.

Non-local Schr\"odinger operators of the type
\begin{equation}
H^\Psi = \Psi(-\Delta) + V
\label{bern}
\end{equation}
currently receive an increasing attention displaying interesting features from the perspective of functional
analysis, stochastic processes and mathematical physics, with a rich interplay between these aspects
\cite{CMS,bib:J,HIL09,KL12a,KL12b,KL15,JW}. Here $\Psi$ is a so called Bernstein function, which can be represented
by the L\'evy-Khintchine formula for subordinators, and $V$ is a Borel-measurable function acting as a
multiplication operator, called potential. For a specific choice of $\Psi$ one obtains fractional Schr\"odinger
operators of the form $(-\Delta + m^{2/\alpha})^{\alpha/2} - m + V$, $0 < \alpha < 2$, with or without rest mass
$m$. A basic motivation of the study of models with fractional Laplacians comes from models of relativistic quantum
mechanics and anomalous transport theory. Another fact that makes these operators particularly interesting is that
they are Markov generators of jump processes.

One approach to investigating the spectral properties of $H^\Psi$ exploits the fact that (\ref{bern}) is a
perturbation of a pseudo-differential operator, and uses direct methods of analysis to derive various spectral
and regularity properties. Another approach, see \cite{HIL09,HL12,KL12a,KL12b,KL15}, uses a probabilistic language
starting from the fact that non-local operators of the type $-\Psi(-\Delta)$ generate Brownian motion sampled at
random times, i.e., subordinate Brownian motions $(B_{S^\Psi_t})_{t\geq 0}$, where the subordinator $\pro {S^\Psi}$
is uniquely determined by the Bernstein function $\Psi$, see \cite{SSV}. Subordinate Brownian motion is a jump
L\'evy process, and when $V$ is added to the operator, it has the effect of killing or reinforcing paths according
to its sign and magnitude. In particular, fractional Schr\"odinger operators with $m=0$ or $m > 0$ generate
rotationally symmetric $\alpha$-stable and relativistic $\alpha$-stable processes, respectively. The two cases
differ essentially by the concentration properties of the L\'evy measure of the underlying process, which has a
significant impact on some properties of eigenfunctions of the related non-local Schr\"odinger operators.

In this paper we take the analytic approach, which complements our previous investigations of non-local Schr\"odinger
operators using a stochastic description. Although the literature in this field is rapidly increasing, there is a
shortage of a detailed understanding of examples which can stand as benchmark cases to the general theory. Our aim
in this paper is to contribute to filling this gap and by choosing $d = 1$, $\alpha = 1$ and $V(x) = x^4$, in the
remainder of this paper we consider the operator (\ref{quart}) in detail. There is special interest in this operator
also from a physics viewpoint in that it describes a semi-relativistic quantum particle with no rest mass, moving in
the force field of a quartic potential. Probabilistically, the kinetic part of $-H$ is the generator of a one-dimensional
Cauchy process (isotropic 1-stable process), which is perturbed by a space-dependent killing by the potential, that can
be derived by the Feynman-Kac formula and subordination.

Using a standard argument \cite[Ch. 3]{LHB} it follows that the spectrum of $H$ is purely discrete consisting of
eigenvalues $0 < \la_1 < \la_2 \leqslant \la_3 \leqslant \ldots \to \infty$, each having finite multiplicity. The
corresponding eigenfunctions $\psi_1, \psi_2, ...$ are bounded continuous functions and form an orthonormal basis
in $L^2(\RR)$. The first eigenfunction (ground state) $\psi_1$, corresponding to the bottom of the spectrum $\lambda_1$,
can be shown to be unique and have a strictly positive version. Making use of recent general results, it follows by
\cite[Cor. 2.2]{KL12b} that for large enough $|x|$ the ground state has a polynomial decay at infinity like
\begin{equation}
\label{-6}
\psi_1(x) \asymp \frac{1}{V(x)|x|^{d+\alpha}} = \frac{1}{|x|^6}.
\end{equation}
Furthermore, by \cite[Cor. 2.1]{KL12b} we have for the other eigenfunctions that
\begin{equation}
|\psi_n(x)| \leqslant C_n \psi_1(x), \quad x \in \RR, \; n \in \NN,
\label{gsdom}
\end{equation}
with a suitable constant $C_n$. (We note that such a pointwise ground state domination fails to hold in general
for usual Schr\"odinger operators, for instance, in the case of $-\Delta + x^2$.) An implication of this is that
the heat kernel $u(t,x,y) = e^{-tH}(x,y)$ asymptotically factorizes and takes the shape of the ground state
exponentially quickly in the sense that
$$
\left|\frac{e^{\lambda_1 t}u(t,x,y)}{\psi_1(x)\psi_1(y)}-1\right| \leqslant C e^{-(\lambda_2-\lambda_1) t},
$$
for sufficiently large $t$ and uniformly in the space variables, at a rate given by the spectral gap $\lambda_2 -
\lambda_1$. In particular, we obtain that $u(t,x,y) \asymp \frac{e^{-\lambda_1t}}{(1+x^6)(1+y^6)}$ for large enough
$t>0$ and all $x,y \in \RR$.

Our goal here is to go beyond these results and derive explicit formulae on the spectrum and eigenfunctions of $H$.
A previous work aiming to solve explicitly the eigenvalue problem for a specific fractional Schr\"odinger operator
has been carried out in \cite{LM12} where the operator $(-\ud^2/\ud x^2)^{1/2} + x^2$, i.e., the case of a quadratic
potential has been considered. In this case the spectrum was found to be the alternating sequence given by the zeroes
of the Airy function $\Ai$ and its derivative $\Ai'$. Also, we have obtained closed formulae for the Fourier transform
of the eigenfunctions in terms of Airy functions. These results allowed to upgrade the general estimates on the fall-off
of eigenfunctions obtained in \cite{KL12b} to tighter bounds and a full asymptotic expansion, and derive detailed
properties of the ground state, heat kernel, heat trace, and spectral gaps.

The features that emerge in the case of the quartic potential differ from the quadratic case on several counts
and we observe a new phenomenon. In the quartic case we can express the Fourier transform of the eigenfunctions
in terms of special functions derived from the fourth-order Airy function of the first kind, denoted
below by $\Ai_4$. Unlike in the quadratic case, where the Fourier transforms are expressed in a single term
involving the Airy function, in the quartic case we have two. One is a highly oscillatory integral given in
terms of $\Ai_4$, while the second, which we will denote by $\widetilde\Ai_4(y)$, is comparable to
$y^{-3/8} e^{-\frac{2\sqrt{2}}{5} y^{5/4}}$ for $y > 0$ and to $|y|^{-3/8} e^{\frac{4}{5} |y|^{5/4}}$ for $y < 0$.
A combination of these functions in the expression of the Fourier transform of $\psi_n$ produces a function dominated
by $\Ai_{4}$ even for low values of $n$, with a subtle contribution from $\widetilde\Ai_4$. The net effect is that
the spectrum $\lambda_n$ of $H$ is located exponentially close with increasing $n$ to the negative zeroes $\mu_n$ of
$\Ai_4$ and $\Ai_4'$ in an alternating order, and $\psi_n$ is exponentially well approximated by the inverse Fourier
transform of $\Ai_4$ at argument shifted by $\la_n$.

We note that numerical and partially rigorous evidence supports our conjecture that this small effect persists for
the higher order anharmonic non-local oscillators  $(-\ud^2/\ud x^2)^{1/2} + x^{2k}$, $k = 3, 4, \ldots$ While even
for the $k=3$ case the expressions rapidly become more complex, the case formally obtained in the limit $k \to \infty$
can be put in a new light through this approach. Note that in this limit $V$ becomes an infinitely deep potential well
with boundaries at $\pm 1$, and the problem becomes equivalent with the non-local Dirichlet problem for
$(-\ud^2/\ud x^2)^{1/2}$ in the interval $(-1,1)$, i.e., the Cauchy process killed outside this interval. It can be
expected that studying the small spectral (non-uniform) shift observed in the present paper will lead to further
understanding of the infinite well problem for which thus far only approximate solutions are around \cite{bib:KKMS}.
Our results in this direction will be further discussed elsewhere.

The remainder of this paper is organized as follows. In Section 2 we first state and then by Fourier transform
reformulate the eigenvalue problem. The so obtained ODE can then be reduced to a fourth-order Airy equation with
specific boundary conditions, and we use suitable modifications of $\Ai_4$ to express the solutions of this
problem (formulae (\ref{eq:fouriereigenfn01c})-(\ref{eq:fouriereigenfn01d})). In Section 3 we identify the spectrum
as the zeroes of higher transcendental functions derived from $\Ai_4$ (Theorem \ref{thm:eigenvalexpressions01}) and
find an expression of these functions in terms of generalized Fresnel integrals (Theorem \ref{thm:eigenvalexpr01}).
Using asymptotic expansions of $\Ai_4$, we derive formulae showing more explicitly the dependence of $\lambda_n$ on
$n$. Next we analyze the small spectral effect discussed above and establish exponential bounds on the differences
between the spectrum and the zeroes of  $\Ai_4$ and $\Ai_4'$ (Theorem \ref{thm:deviation01}). We also obtain bounds
on the full sequence of spectral gaps $\lambda_{n+1} - \lambda_n$. In a last subsection we derive the behaviour of
the heat trace around the origin, and obtain a Weyl-type law on the asymptotic behaviour of the counting measure of
the eigenvalues (Theorem \ref{thm:traceest}, Corollary \ref{thm:weylasymp}). In Section 4 we turn to discussing the
eigenfunctions. Since fourth-order Airy functions appear to be little explored in the literature, we start by a
brief presentation of their properties particularizing the results obtained in \cite{D14}. To keep the length of
this paper reasonable while offering detail, we only discuss the material which will be strictly needed here and
refer for more to the given source. Next we obtain expressions of the Fourier transforms of the eigenfunctions
in terms of these special functions (Theorem \ref{thm:eigenfns01}), and also prove analyticity of the eigenfunctions
(Theorem \ref{thm:analytic01}). In Theorem \ref{thm:diffsupbound} we analyze the small effect on the level of
eigenfunctions and show that the $\psi_n$ differ by exponentially small errors from the inverse Fourier transform
of the $L^2$-normalized $\Ai_4$ at values shifted by $\lambda_n$. Proofs will be presented in Section 5.
%In separate sections of the Appendix we show tables and plots illustrating the relevant special functions, spectrum and eigenfunctions.

\section{Eigenvalue problem for the quartic oscillator}\label{sec:spectralanal01}
Let $H_0 = (-\ud^2/\ud x^2)^{1/2}$ be the square root of the Laplace operator in one dimension, defined by
$$
({\mathcal F}(-\ud^2/\ud x^2)^{1/2}f)(y) = |y|({\mathcal F} f)(y)
$$
with domain $H^1(\RR)= \{ f \in L^2(\RR): |y|{\mathcal F} f \in L^2(\RR)\}$, where $\mathcal F$ denotes Fourier
transform. $H_0$ is essentially self-adjoint on $C_0^\infty(\RR)$ and $\Spec H_0 =\Spec_{\rm{ess}} H_0 =[0,\infty)$.

In this paper we consider the fractional Schr\"odinger operator formally written as $H = H _0 + x^4$. One way to
define this operator as a self-adjoint operator with dense domain in $L^2(\RR)$ is by using the Feynman-Kac formula
$$
(f, e^{-tH} g)_{L^2(\mathbb R)} = \int_{-\infty}^\infty f(x) \ex^x[e^{-\int_0^t X_s^4 ds} g(X_t)] dx, \quad t \geq 0,
\, f, g \in L^2(\RR),
$$
where $\pro X$ is a one-dimensional Cauchy process, and the expectation at the right hand side is taken with respect
to the measure of this process starting from $x \in \RR^d$. The right hand side can be proven to be a strongly continuous
semigroup, and thus it follows by the Hille-Yoshida theorem that $H$ is self-adjoint and equals the form sum $H_0
\dot {+} x^4$. Since our framework here is other than using functional representation, we refer for further details
to \cite[Th. 4.8]{HIL09}. As the potential is positive and growing to infinity, we have as a consequence of Rellich's
theorem \cite[Th. 3.20]{LHB} that $H$ has compact resolvent, in particular, $\Spec H \subset [0,\infty)$ is discrete,
consisting of isolated eigenvalues, each of finite multiplicity. Our concern in what follows is to determine and analyze
this spectrum and the corresponding eigenfunctions.

Consider the eigenvalue equation
\eq{eq:quarticoper01a}
 \left(-\frac{\ud^{2}}{\ud x^{2}}\right)^{1/2}\psi_n + x^{4}\psi_n = \lambda_n \psi_n,
\quad \psi_n \in H^1(\RR), \; n \in \Bbb N.
\en
From relations (\ref{-6})-(\ref{gsdom}) it directly follows that
\begin{equation}
\label{L1}
x^4 \psi_n \in L^1(\RR), \quad n \in \NN.
\end{equation}
This implies, in particular, $\psi_n \in L^1(\RR)$ for every $n \in \NN$. We write $\phi_n(y) := ({\mathcal F}\psi_n)(y) =
\frac{1}{\sqrt{2\pi}}\int_{-\infty}^{\infty} e^{- i x y } \psi_n(x) \ud x$, $y \in \RR$, for the Fourier transform of
$\psi_n$. Since
$$
\left({\mathcal F} \lk -\frac{\ud^{2}}{\ud x^{2}}\rk^{1/2} {\mathcal F}^{-1}\right)(y) = |y| \quad \mbox{and} \quad
({\mathcal F} x^4 {\mathcal F}^{-1})(y) = \frac{\ud^4}{\ud y^4},
$$
equation \eqref{eq:quarticoper01a} transforms into
\eq{eq:diff01a}
\frac{\ud^4}{\ud y^4}\phi_n(y)+(|y|-\la_n) \phi_n(y)=0, \quad y \in \RR.
\en
By (\ref{L1}) we have
\begin{align}\label{eq:eigenfnfourier01}
\Big|\frac{\ud^{4}}{\ud y^{4}}\phi_{n}(y)\Big| &=  \frac{1}{\sqrt{2\pi}}\bigg|\int_{-\infty}^{\infty} x^{4} e^{- i x y }
\psi_{n}(x) \ud x\bigg| \leqslant \frac{1}{\sqrt{2\pi}} \|x^4\psi_{n}\|_1 < \infty, \quad y \in \RR.
\non
\end{align}
Furthermore, by multiplying and dividing, and using the Schwarz inequality
\begin{align*}
\int_{-\infty}^\infty |\phi_n(y)| dy \leqslant
\left(\int_{-\infty}^\infty \frac{dy}{1+y^2}\int_{-\infty}^\infty (1+y^2)|\phi_n(y)|^2 dy \right)^{1/2}=
\sqrt\pi \Vert \psi_{n} \Vert_{H^{1}(\RR)}.
\end{align*}
From these two observations we have that $\phi_n \in C^4(\RR) \cap L^1(\RR)$, for all $n\in \NN$.

Note, moreover, that if $\phi_{n}(y)$ is a solution of \eqref{eq:diff01a}, then so is $\phi_{n}(-y)$.
Thus it suffices to consider \eqref{eq:diff01a} only for $y > 0$, and construct odd and even solutions on
the whole of $\RR$. They respectively satisfy
\eq{eq:bdrycond01}
\phi_{n}(0+) = \phi_{n}^{\prime \prime}(0+) =0, \quad n = 2, 4, 6, \ldots
\en
and
\eq{eq:bdrycond02}
\phi_{n}^{\prime}(0+) = \phi_{n}^{\prime \prime \prime}(0+) =0, \quad n = 1, 3, 5, \ldots
\en
Our goal translates then to studying the spectrum and the eigenfunctions of $H$ given by \eqref{eq:quarticoper01a}
through the $L^1$ solutions of equation \eqref{eq:diff01a} satisfying the boundary conditions
\eqref{eq:bdrycond01}-\eqref{eq:bdrycond02}.

Let $y > 0$. By the change of variable $y-\la_{n} \mapsto y$ and using the notation $\varphi_n(y ) \equiv \phi_{n}(y+\la_{n})$,
we obtain
\eq{eq:diff01}
\frac{\ud^4}{\ud y^4}\varphi_n(y)+y \varphi_n(y)=0.
\en
Equation \eqref{eq:diff01} can be thought of being a higher order Airy equation. Following the terminology in
\cite{bib:koh79} we call the particular solution
\eq{eq:extendAiry04}
\Ai_{4}(y) := \frac{1}{\pi}\int_{0}^{\infty} \cos\lk \frac{t^5}{5} + y t \rk \ud t
\en
a \emph{fourth-order Airy function}. These functions seem to be relatively little known in the literature and we will
use Section \ref{sec:extairy} to presenting some of their basic properties used in this paper, see \cite{D14} for
details and a more general study.

Notice that equation \eqref{eq:diff01} is invariant under the rotations $y \mapsto y e^{2l\pi i/5}$, $l \in \ZZ$, thus
we have four other solutions by rotation of the argument of $\Ai_{4}(y)$, i.e.,
\eq{eq:soln01}
\Ai_{4}(y e^{-2\pi i/5}), \quad \Ai_{4}(y e^{-4\pi i/5}), \quad  \Ai_{4}(y e^{4\pi i/5}), \quad \Ai_{4}(y e^{2\pi i/5}).
\en
Any four of the above five solutions can be shown to have a non-vanishing Wronskian.
%, thus a particular solution of \eqref{eq:diff01} can be constructed by taking linear combinations.
Consider
\begin{eqnarray*}
&&
\rA_{1}(y) := e^{-2\pi i/5} \Ai_{4}(y e^{-2\pi i/5}), \quad \rA_{2}(y) := e^{-4\pi i/5} \Ai_{4}(y e^{-4\pi i/5}) \\
&&
\rA_{3}(y) := e^{4\pi i/5} \Ai_{4}(y e^{4\pi i/5}), \quad \rA_{4}(y) := e^{2\pi i/5} \Ai_{4}(y e^{2\pi i/5}).
\end{eqnarray*}
Cauchy's theorem on contour integration then gives
$\Ai_{4}(y) = -\sum_{r=1}^{4} \rA_{r}(y)$.
%Three other solutions independent of $\Ai_{4}(y)$ can also be constructed by taking other suitable linear combinations
%of $\rA_{r}(y)$ for $r=1,2,3,4$, which together with $\Ai_{4}(y)$ form a complete set of independent solutions to the
%differential equation \eqref{eq:diff01}. They are given by $\widetilde{\Ai}_{4}(y) = i\sum_{r=1}^{4} (-1)^{r+1} \rA_{r}(y)$
%with integral representation (see Appendix \ref{sols} for details)
By taking $\Ai_{4}(y)$ and suitable linear combinations of $\rA_{r}(y)$ for $r=1,2,3,4$, a complete set of independent
solutions to the differential equation \eqref{eq:diff01} is given by $\widetilde{\Ai}_{4}(y) = i\sum_{r=1}^{4} (-1)^{r+1}
\rA_{r}(y)$ with integral representation %(see Appendix \ref{sols} for details)
\eq{eq:secondef}
\widetilde{\Ai}_{4}(y) = \frac{1}{\pi}\int_{0}^{\infty} \lkkk e^{-yt - \frac{t^{5}}{5}} - \sin\lk \frac{t^{5}}{5} + yt
\rk \rkkk \ud t,
\en
and ${\rm G}_{3}(y) =  \rA_{1}(y) + \rA_{4}(y) - \lk \rA_{2}(y) + \rA_{3}(y) \rk$, ${\rm G}_{4}(y) = i\left(\rA_{4}(y)
- \rA_{1}(y)\right) + i\left(\rA_{3}(y) - \rA_{2}(y) \right)$. Since ${\rm G}_{3}, {\rm G}_{4}$ diverge as $y \to
\infty$, these solutions are ruled out by the $L^{1}$ condition. Hence finally we obtain that the even and odd solutions
of equation \eqref{eq:diff01a} satisfying \eqref{eq:bdrycond01}--\eqref{eq:bdrycond02} are
\begin{eqnarray}
\label{eq:fouriereigenfn01c}
&&\phi_{2j-1}(y) =  c_{1,2j-1} \Ai_{4}(|y|-\la_{2j-1}) + c_{2,2j-1} \widetilde{\Ai}_{4}(|y|-\la_{2j-1})\\
\label{eq:fouriereigenfn01d}
&& \phi_{2j}(y) = c_{1,2j} {\rm sgn}(y)\Ai_{4}(|y|-\la_{2j}) + c_{2,2j} {\rm sgn}(y)\widetilde{\Ai}_{4}(|y|-\la_{2j}),
\end{eqnarray}
for all $y \in \RR$, $j \in \Bbb N$, with constant prefactors to be determined below. The above calculations are simple
and we leave the details to the reader; alternatively see \cite{D14}.

\section{The spectrum of $H$}\label{sec:spectralanal02}
\subsection{Identification of the spectrum}\label{sec:spectrumid}
The boundary conditions \eqref{eq:bdrycond01}--\eqref{eq:bdrycond02} applied to \eqref{eq:fouriereigenfn01c}
and \eqref{eq:fouriereigenfn01d} yield for $n = 2, 4, 6, \ldots$ and $n = 1, 3, 5, \ldots$, respectively, that
\eq{eq:det01}
{\rm det} \begin{pmatrix}
\Ai(-\la_{n}) & \widetilde{\Ai}_{4}(-\la_{n}) \\
\Ai^{\prime \prime}(-\la_{n}) & \widetilde{\Ai}_{4}^{\prime \prime}(-\la_{n})
\end{pmatrix} = 0
%\en
%for $n = 2, 4, 6, \cdots$, and
%\eq{eq:det02}
\qquad \mbox{and} \qquad
{\rm det}\begin{pmatrix}
\Ai^{\prime}(-\la_{n}) & \widetilde{\Ai}_{4}^{\prime}(-\la_{n}) \\
\Ai^{\prime \prime \prime}(-\la_{n}) & \widetilde{\Ai}_{4}^{\prime \prime \prime}(-\la_{n})
\end{pmatrix} = 0.
\en
Define
\eq{eq:eigenvalexpr01a}
\Phi_{1}(y) := \Ai(y)\widetilde{\Ai}_{4}''(y) - \Ai''(y)\widetilde{\Ai}_{4}(y) \quad \mbox{and} \quad
\Phi_{2}(y) := \Ai'(y)\widetilde{\Ai}_{4}'''(y) - \Ai'''(y)\widetilde{\Ai}_{4}'(y).
\en
Then we have the following main result of this section.
\bt{thm:eigenvalexpressions01}
%The spectrum of $H$ is purely discrete with isolated eigenvalues given by
\begin{align}
\label{ev}
\Spec H = \{\lambda_n, \, n \in \mathbb N: \, \la_{2j-1} = -a_{2,j} \;\, {\rm and} \;\, \la_{2j} = -a_{1,j}, \,
j \in \Bbb N\},
\end{align}
where $a_{1,j}$ and $a_{2,j}$ are the negative real zeroes of the higher transcendental
functions $\Phi_{1}$ and $\Phi_{2}$, respectively, arranged in increasing order.
\et

The following result allows to identify the eigenvalues as the zeroes of more familiar special functions.
Recall the generalised Fresnel sine and cosine integrals ${\rm si}(a,z)$ and ${\rm ci}(a,z)$ (see also
(\ref{fresnel}) below). The proof of this and the forthcoming results will be presented in Section 5.
\bp{thm:eigenvalexpr01}
For all $\lambda > 0$, we have the expressions
\eq{eq:eigenvalexpr01b}
\Phi_{1}(-\la) = -\frac{1}{2\pi^{2}} \int_{0}^{\infty} \lkkk h_{1}(v) \cos\lk \frac{v^{5}}{20} +
\la v \rk -  h_{2}(v)\sin\lk \frac{v^{5}}{20} + \la v \rk \rkkk \ud v
\en
and
\eq{eq:eigenvalexpr02b}
\Phi_{2}(-\la) = -\frac{1}{2\pi^2} \int_{0}^{\infty} \lkkk h_{3}(v) \cos\lk \frac{v^{5}}{20} +
\la v \rk - h_{4}(v) \sin\lk \frac{v^{5}}{20} + \la v \rk \rkkk \ud v,
\en
where
\begin{align*}
h_{1}(v) &:= v^{1/2}{\rm si}\lk \frac{1}{2}, \frac{v^{5}}{16} \rk, \quad
h_{2}(v) := v^{1/2}{\rm ci}\lk \frac{1}{2}, \frac{v^{5}}{16} \rk \\
h_{3}(v) &:= 2\sin\frac{v^5}{16} - v^{5/2} {\rm ci}\lk \frac{1}{2}, \frac{v^5}{16} \rk, \quad
h_{4}(v) := 2\cos\frac{v^5}{16} + v^{5/2} {\rm si}\lk \frac{1}{2}, \frac{v^5}{16} \rk.
\end{align*}
\ep

\medskip
\subsection{Approximations of the spectrum}\label{sec:specapprox}
To have a more specific idea of the dependence of the eigenvalues on $n$, we can use Proposition \ref{thm:eigenvalexpr01}
to derive asymptotic relations for $\Phi_{1}(-\la)$ and $\Phi_{2}(-\la)$.
\bc{cor:asymptrelation}
There exists $0 < \la_{0} \leqslant 1/2$ such that for every $\la > \la_{0}$, we have
\eq{eq:asymprel01}
\Phi_{1}(-\la) = \frac{1}{\pi}\la^{-1/4} e^{\frac{4}{5}\la^{5/4}}\sin\lk \frac{4}{5}\la^{5/4} +
\frac{\pi}{4} \rk + O\lk \frac{e^{\frac{4}{5}\la^{5/4}}}{\la^{7/8} } \rk
\en
and
\eq{eq:asymprel02}
\Phi_{2}(-\la) = \frac{1}{\pi}\la^{1/4} e^{\frac{4}{5}\la^{5/4}}\cos\lk \frac{4}{5}\la^{5/4} +
\frac{\pi}{4} \rk + O\lk \frac{e^{\frac{4}{5}\la^{5/4}}}{\la^{3/8} } \rk.
\en
\ec

Solving the equations $\Phi_{1}(-\la_{n}) = 0$ for $n = 2j$ and $\Phi_{2}(-\la_{n}) = 0$ for $n = 2j-1$,
$j \in \Bbb N$, the above asymptotic relations give some more explicit expressions for the eigenvalues.
\bc{cor:eigenvalexpressions02}
For each $n \in \NN$,
\eq{eq:eigenvalrel02}
\la_{n} = \lk \frac{5(2n-1)\pi}{16} \rk^{4/5} \lkkk 1+O\lk \frac{1}{n^{3/2}} \rk \rkkk.
\en
\ec

The next result describes the small spectral effect discussed in the Introduction, i.e., that the eigenvalues $\la_{n}$
differ from the zeroes of $\Ai_{4}$ and $\Ai_4'$ by exponentially small terms. Define
\eq{eq:zerosrep}
\mu_{n} :=
\begin{cases}
-\ap^{\prime}_j \quad & n = 2j-1 \\
-\ap_j \quad & n = 2j
\end{cases}
\en
for $j \in \NN$, where $\ap_j$ and $\ap_j'$ are the negative real zeroes of $\Ai_4$ and $\Ai_4'$, respectively, see
\eqref{eq:extairyzer01}-\eqref{eq:extairyzer01b} below.
\bt{thm:deviation01}
For every $n \in \mathbb N$ we have that
\eq{eq:devineq01}
|\la_{n} - \mu_{n}| \leqslant \frac{2^{5/4}\sqrt{\GA\lk \frac{4}{5}\rk}}{5^{1/10}\sqrt{\pi}} \, \ga_{n}
\en
where $\ga_{n} := \sqrt{c_{1,n} c_{2,n}}$ and $c_{1,n}, c_{2,n}$ are the normalization constants given in Theorem
\ref{thm:eigenfns01} below. Moreover,
$$
\ga_{n} \asymp \sqrt{\frac{\pi}{2}}\la_{n}^{-1/8}e^{-\frac{2}{5}\la_{n}^{5/4}} \sim \sqrt{\frac{\pi}{2}} \,
n^{-1/10}e^{-\frac{2}{5}n,}
$$
for all $n \in \mathbb N$.
\et
\noindent
We note that throughout this paper the standard Landau notations are used, and $f(x) \sim g(x)$ means $f(x) = g(x)[1+o(1)]$,
while $f(x) \asymp g(x)$ means that there exist real numbers $C_{2} \geqslant C_{1} > 0$ such that $C_{1} g(x) \leqslant f(x)
\leqslant C_{2} g(x)$.

This gives the following bounds on the sequence  $\la_{n+1} - \la_{n}$ of spectral gaps.
\bc{thm:specgaps01}
For every $n \in \NN$ we have that
\eq{eq:specgaps01b}
\frac{\pi}{2}\lk \frac{8}{15\pi} \rk^{1/5} \lk n - \frac{1}{2} \rk^{-1/5} \leqslant \la_{n+1} - \la_{n} \leqslant
\frac{\pi}{2}\lk \frac{8}{5\pi} \rk^{1/5} \lk n - \frac{1}{2} \rk^{-1/5}.
\en
\ec
\noindent

\medskip
\subsection{Heat trace and Weyl-type theorem}\label{subsec:spectralanal03}
Using the asymptotic expression in Corollary \ref{cor:eigenvalexpressions02}, we are able to derive an expression
for the behaviour of the trace of the semigroup
$$
Z(t):= \Tr e^{-tH} = \sum_{n=1}^{\infty} e^{-\la_{n}t}, \quad t > 0,
$$
in a neighbourhood of the origin.
\bt{thm:traceest}
We have
\eq{eq:traceasympformula01}
\lim_{t \downarrow 0}t^{5/4}Z(t) = \frac{2\GA\lk \frac{5}{4}\rk}{\pi}.
\en
\et
A consequence of this is the following Weyl-type asymptotic formula on the distribution of eigenvalues.
Denote by
$$
N(\la) := |\{n \in \mathbb N: \lambda_n \leq \lambda\}|
$$
the counting measure of the number of eigenvalues $\la_{n}$ of $H$ up to level $\la > 0$.

\begin{corollary}
\label{thm:weylasymp}
We have that
\eq{eq:weylasyformula01}
\lim_{\la \to \infty} \frac{N(\la)}{\la^{5/4}} = \frac{8}{5\pi}.
\en
\end{corollary}

\bigskip

\section{Expressions and properties of eigenfunctions}\label{subsec:spectralanal04}
\subsection{Fourth-order Airy function}\label{sec:extairy}
The fourth-order Airy function of the first kind $\Ai_{4}$ is a particular case of the more general class $\Aik$, $k = 2, 4, 6,
\ldots$, formally defined by the integrals
\[
\Aik(y) = \frac{1}{2\pi i} \int_{\cG} e^{-\ga_{k}ys - \frac{s^{k+1}}{k+1}} \ud s, \qquad y \in \RR, \;\; \ga_{k}:= (-1)^{k/2},
\]
where $\cG$ is any infinite contour in the complex $s$-plane that starts at infinity in the sector $-\frac{k\pi}{k+1} - \del <
{\rm arg}(s) < \del - \frac{k\pi}{k+1}$ and ends at infinity in the sector $\frac{k\pi}{k+1} - \del < {\rm arg}(s) <
\frac{k\pi}{k+1} + \del$ for $0 \leqslant \del < \frac{(k-1)\pi}{2(k+1)}$, cutting through the negative real semi-axis. The integral
is convergent since $\Re e(s^{k+1}) \to \infty$ as $|s| \to \infty$ within $\big|{\rm arg}(s) \pm \frac{k\pi}{k+1}\big|< \del$. It
is also possible to continuously deform the contour $\cG$ to align with the imaginary axis without altering the value of the integral,
and hence by the change of variable $s=it$ and for $k = 4$, we obtain \eqref{eq:extendAiry04}. The function $\Aik$ is smooth and
bounded, and extends to the complex plane as an entire function, and hence it follows directly that $\Aik$ satisfies the higher
order Airy differential equation
\[
\frac{\ud^{k}}{\ud y^{k}} \varphi(y) + \ga_{k} y \varphi(y) = 0.
\]
Since proofs involve lengthy calculations, for more details we refer to  \cite{D14}; here we present some properties of
$\Ai_{4}$ taken from here for the specific case $k = 4$ which are of particular importance in this paper.

\bp{thm:analexpansionzero}
The following expansion holds:
\eq{airyseriesexpansion}
\Ai_{4}(y) = \sum_{p=0}^{3} \frac{\Ai_{4}^{(p)}(0)}{p!} y^{p} \sum_{l = 0}^{\infty} \frac{(-1)^{l}}{l!}
\prod_{\substack{j=0 \\ j \neq p}}^{3} \frac{\GA\lk \frac{p-j}{5} +1\rk}{\GA\lk \frac{p-j}{5} + l +1\rk} \frac{y^{5l}}{5^{4l}},
\quad  y \in \RR.
\en
Moreover, for the derivatives we have
\eq{eq:initialvalues}
\Ai_{4}^{(p)}(0) = \frac{\cos\lk \frac{\pi}{2}  \frac{p+1}{5} + \frac{p\pi}{2} \rk}{5^{\frac{4-p}{5}}\GA\lk \frac{4-p}{5} \rk\sin\lk
\frac{p+1}{5}\pi \rk}, \quad \mbox{for each $p=0,1,2,3$}.
\en
\ep
\noindent
Formula \eqref{airyseriesexpansion} can alternatively be expressed as the sum
\eq{eq:sumshypergeo01}
\Ai_{4}(y) = \sum_{p=0}^{3} \Ai_{4}^{(p)}(0) \frac{y^{p}}{p!} {\rm U}_{p}(y),
\en
in terms of the generalized hypergeometric functions $_0\rF_{3}(\ap, \BT, \ga; z)$, where we denote
${\rm U}_{0}(y) :=  \,_0\rF_{3}\lk \frac{2}{5},\frac{3}{5},\frac{4}{5};-\frac{y^5}{5^4}\rk$,
${\rm U}_{1}(y) :=  \,_0\rF_{3}\lk \frac{3}{5},\frac{4}{5},\frac{6}{5};-\frac{y^5}{5^4}\rk$,
${\rm U}_{2}(y) :=  \,_0\rF_{3}\lk \frac{4}{5},\frac{6}{5},\frac{7}{5};-\frac{y^5}{5^4}\rk$, and
${\rm U}_{3}(y) :=  \,_0\rF_{3}\lk \frac{6}{5},\frac{7}{5},\frac{8}{5};-\frac{y^5}{5^4}\rk$.

\bp{prop:tailbehav}
There exist constants $C_{1}, C_{2} > 0$ such that
\[
|\Ai_{4}(y)| \leqslant \min \lkk C_{1} e^{-2y}, C_{2}|y|^{-\frac{5}{8}}\rkk, \quad y \in \RR.
\]
Moreover, we have that $\Ai_{4}\rceil_{\RR^+} \in L^{1}(\RR^+)\cap L^{2}(\RR^+)$.
\ep

The next result gives the behaviour of $\Ai_{4}$ and $\Ai^{\prime}_{4}$ on the negative semi-axis from which
the asymptotic expansions of the negative real zeroes of $\Ai_{4}(y)$ and $\Ai^{\prime}_{4}(y)$ are derived.
\bp{thm:airyasymp02}
We have that
$$
\Ai_{4}(-y) = \frac{1}{\pi^{1/2} y^{3/8}}\lk \cos\lk \xi - \frac{\pi}{4} \rk P\lk \xi \rk +
\sin\lk \xi - \frac{\pi}{4} \rk Q\lk \xi \rk \rk
$$
and
$$
\Ai^{\prime}_{4}(-y) = \frac{1}{\pi^{1/2}y^{1/8}}\lk \sin\lk \xi - \frac{\pi}{4} \rk R\lk \xi \rk -
\cos\lk \xi - \frac{\pi}{4} \rk S\lk \xi \rk \rk,
$$
where
$$
P\lk \xi \rk \sim \sum_{r=0}^{\infty}(-1)^{r}\frac{C_{2r}}{\xi^{2r}} \quad \mbox{and} \quad
Q\lk \xi\rk \sim \sum_{r=0}^{\infty}(-1)^{r}\frac{C_{2r+1}}{\xi^{2r+1}} \quad \mbox{as $\xi \to +\infty$}
$$
and
$$
R\lk \xi \rk \sim \sum_{r=0}^{\infty}(-1)^{r}\frac{C^{\prime}_{2r}}{\xi^{2r}} \quad \mbox{and} \quad S\lk \xi\rk
\sim \sum_{r=0}^{\infty}(-1)^{r}\frac{C^{\prime}_{2r+1}}{\xi^{2r+1}} \quad \mbox{as $\xi \to +\infty$},
$$
with $\xi = \frac{4}{5}y^{5/4}$, and coefficients $C_{l}, C^{\prime}_{l}$ expressed in terms of the Bell polynomials.
\ep
Finally, we give some asymptotic expansions of the negative real zeroes of $\Ai_{4}$ and $\Ai_{4}^{\prime}$ denoted by
$\ap_{n}$ and $\ap_{n}^{\prime}$, respectively.
\bp{thm:zeros}
For every $n \in \NN$ we have
\eq{eq:extairyzer01}
\ap_{n} = -f\lk \frac{5(4n - 1)\pi}{16} \rk  \quad \mbox{and} \quad \ap^{\prime}_{n} = -g\lk \frac{5(4n - 3)\pi}{16} \rk,
\en
where
\eq{eq:extairyzer01a}
f(\tau) \sim \tau^{4/5}\lk 1 + \frac{0.1586783204}{\tau^2} - \frac{0.03595263992}{\tau^4} + \frac{0.01511323043}{\tau^6}
-  \cdots \rk ,
\en
\eq{eq:extairyzer01b}
g(\tau) \sim \tau^{4/5}\lk 1 - \frac{0.05289277344}{\tau^2} - \frac{0.003797437462}{\tau^4} - \frac{0.0005042991615}{\tau^6} - \cdots \rk,
\en
as $\tau \to \infty$.
\ep

\medskip
\subsection{Eigenfunctions}
The Fourier transforms of the eigenfunctions of $H$ can be expressed in terms of fourth-order extended Airy functions. Define
\begin{gather}
\LA_{1}(\xi) := \widetilde{\Ai}_{4}(\xi)\Ai_{4}^{\prime}(\xi) - \Ai_{4}(\xi) \widetilde{\Ai}_{4}^{\prime}(\xi)
\label{eq:constdem01}\\
\LA_{2}(\xi) := \widetilde{\Ai}_{4}^{\prime}(\xi)\Ai_{4}^{\prime \prime}(\xi) -
\Ai_{4}^{\prime}(\xi) \widetilde{\Ai}_{4}^{\prime \prime}(\xi)
\label{eq:constdem02} \\
\LA_{3}(\xi) := \widetilde{\Ai}_{4}^{\prime \prime}(\xi) \Ai_{4}^{\prime \prime \prime}(\xi) -
\Ai_{4}^{\prime \prime}(\xi) \widetilde{\Ai}_{4}^{\prime \prime \prime}(\xi).
\label{eq:constdem03}
\end{gather}
\bt{thm:eigenfns01}
The Fourier transform of the $L^{2}$-normalized eigenfunctions of $H$ is given by
\eq{eq:fouriereigenfns}
\lk {\mathcal F} \psi_{n} \rk(y) =
\begin{cases}
c_{1,n} \Ai_{4}(|y|-\la_{n}) + c_{2,n}\widetilde{\rm Ai_{4}}(|y|-\la_{n}), \quad &n=2j-1 \\ \\
c_{1,n} {\rm sgn}(y)\Ai_{4}(|y|-\la_{n}) + c_{2,n}{\rm sgn}(y)\widetilde{\rm Ai_{4}}(|y|-\la_{n}), \quad &n=2j
\end{cases}
\en
for all $y \in \RR$ and $j \in \mathbb N$, where
\eq{eq:normconst01}
c_{1,n} :=
\begin{cases}
-\frac{1}{\sqrt{2}} \frac{\widetilde{\Ai}_{4}^{\prime}(-\la_{n})}{\sqrt{\la_{n}\LA_{1}^{2}(-\la_{n}) + \LA_{2}^{2}(-\la_{n})}}
\quad & n = 2j-1 \\
\frac{1}{2} \frac{\widetilde{\Ai}_{4}^{\prime \prime}(-\la_{n})}{\sqrt{\LA_{2}(-\la_{n}) \LA_{3}(-\la_{n})}} \quad & n = 2j
\end{cases}
\en
and
\eq{eq:normconst02}
c_{2,n} :=
\begin{cases}
\frac{1}{\sqrt{2}} \frac{\Ai_{4}^{\prime}(-\la_{n})}{\sqrt{\la_{n}\LA_{1}^{2}(-\la_{n}) + \LA_{2}^{2}(-\la_{n})}}
\quad & n = 2j-1 \\
-\frac{1}{2} \frac{\Ai_{4}^{\prime \prime}(-\la_{n})}{\sqrt{\LA_{2}(-\la_{n}) \LA_{3}(-\la_{n})}} \quad & n = 2j
\end{cases}
\en
\et

We now proceed with further properties of the eigenfunctions $\psi_{n}$. First we obtain a full asymptotic expansion.
\bt{thm:eigenfnasymptotics}
For every $j \in \NN$ and $N=2,3,\ldots$
\begin{eqnarray}
\psi_{2j-1}(x) &=& \sum_{l = 1}^{N-1} (-1)^{l+1} \frac{\cP(\la_{2j-1})}{x^{4+2l}}+O\left(\frac{1}{x^{2N+4}}\right)
\label{eq:eigenfnasymp01a} \\
\psi_{2j}(x) &=& \sum_{l = 1}^{N-1} (-1)^{l+1} \frac{\cQ(\la_{2j})}{x^{5+2l}}+O\left(\frac{1}{x^{2N+5}}\right)
\label{eq:eigenfnasymp01b}
\end{eqnarray}
holds as $x\to\infty$, where we have that $\cP(-\la_{1}) = \sqrt{\frac{2}{\pi}}\left(c_{1,1}\Ai_{4}(-\la_{1}) +
c_{2,1}\widetilde{\Ai}_{4}(-\la_{1})\right)$, $\cQ(-\la_{2}) = 2\sqrt{\frac{2}{\pi}} \left(c_{1,2}\Aip_{4}(-\la_{2}) +
c_{2,2}\widetilde{\Ai}_{4}^{\prime}(-\la_{2})\right)$, and for all $j = 2, 3, \ldots$
\begin{gather*}
\cP(\la_{2j-1}) = \sqrt{\frac{2}{\pi}}\left(c_{1,2j-1}\Ai_{4}^{(3+2l)}(-\la_{2j-1}) + c_{2,2j-1}\widetilde{\Ai}_{4}^{(3+2l)}(-\la_{2j-1})\right) \\
\cQ(\la_{2j}) = \sqrt{\frac{2}{\pi}} \left(c_{1,2j-1}\Ai_{4}^{(4+2l)}(-\la_{2j}) + c_{2,2j}\widetilde{\Ai}_{4}^{(4+2l)}(-\la_{2j})\right).
\end{gather*}
\et
\noindent
%In Appendix \ref{numeref} we present some plots of the first few eigenfunctions obtained numerically.

From general results on the smoothing property of the evolution semigroup $e^{-tH}$ we know that the eigenfunctions are bounded and
continuous functions \cite{KL12b}. Using the extra information for the specific case here, we obtain stronger regularity properties.
\bt{thm:analytic01}
For every $n \in \NN$ the eigenfunction $\psi_{n}$ is analytic on $\RR$ and has the expansion
\begin{equation}
\psi_{n}(x) =
\begin{cases}
\sum_{r=0}^{\infty}(-1)^{r} a_{2r}(\la_{n})x^{2r},  \quad &n=2j-1  \\ \\
\sum_{r=0}^{\infty}(-1)^{r} a_{2r+1}(\la_{n}) x^{2r+1}, \quad & n = 2j
\end{cases}
\end{equation}
for all $j \in \mathbb N$ and $x \in \RR$, where
$$
a_{p}(u) = \frac{1}{p!} \sqrt{\frac{2}{\pi}} \int_{0}^{\infty}y^{p}\left( c_{1,n}\Ai_{4}(y-u) +
c_{2,n}\widetilde{\Ai}_{4}(y-u)\right) \ud y, \quad p \in \mathbb N.
$$
\et

The following result shows another aspect of the small effect discussed before on the level of the spectrum. Its occurrence
on the level of eigenfunctions consists in the fact that the eigenfunctions $\psi_n$ are exponentially well approximated by
the inverse Fourier transform of the dominating term in (\ref{eq:fouriereigenfn01c})-(\ref{eq:fouriereigenfn01d}). Write
$\chi_{n}(y) := c_{1,n} \Ai_{4}(|y| - \mu_{n})$, for which we have $\|\chi_{n}\|_{2}=1$ for all $n \in \NN$. Using that
$|\psi_{n} - \cF^{-1}\chi_{n}| \leqslant \frac{1}{\sqrt{2\pi}}\int_{-\infty}^\infty |\phi_{n}(y) - \chi_{n}(y)|\ud y$,
we derive a direct estimate on the sup-norm of the difference, and another using the local cancellations between successive
eigenvalues.
\bt{thm:diffsupbound}
For every $n \in \mathbb N$ we have
\eq{eq:uniformbdd01}
\Big\Vert \psi_{n} - \cF^{-1}\chi_{n} \Big\Vert_{\infty} \leqslant \frac{2^{9/4}}{\sqrt{5 \pi}} c_{2,n},
\en
where $c_{2,n} \asymp \sqrt{\frac{\pi}{2}} \, \la_{n}^{-1/8} e^{-\frac{4}{5} \la_{n}^{5/4}} \sim \sqrt{\frac{\pi}{2}} \,
n^{-1/10} e^{-\frac{4}{5} n}$. Moreover, for every $n, p \in \NN$ it follows that
\eq{eq:uniformbdd01a}
\lk \int_{-\infty}^{\la_{1}} + \sum_{p \geqslant 1} \int_{\la_{p}}^{\la_{p+1}} \rk \Big|\phi_{n}(y) - \chi_{n}(y) \Big| \ud y
\leqslant \frac{2^{3/4}}{\sqrt{5}} J(n)
\en
where
$$
J(n) =  3 c_{2,n} + B_{n} + \frac{8C_{n}}{5^{5/4} 3^{1/4}} + \lk \frac{1}{e^{\frac{2\pi}{5^{5/4} 3^{1/4}}} - 1}
+ \frac{e^{-\frac{\pi}{15^{1/4} 2}}}{e^{\frac{\pi}{15^{1/4} 2}} - 1} \rk D_{n} - \frac{e^{\frac{\pi}{15^{1/4} 2}} E_{n}}
{\lk e^{\frac{\pi}{15^{1/4} 2}} - 1\rk}
$$
with
\begin{gather*}
B_{n} = e^{-\frac{2^{7/4}\pi \, n}{5^{5/4} 3^{1/4}}} c_{2,n} \sim \sqrt{\frac{\pi}{2}} \, n^{-1/10} e^{-1.8739 \,n}  \quad
\mbox{and} \quad
C_{n} = n^{-1/4} c_{2,n} \sim \sqrt{\frac{\pi}{2}} \, n^{-7/20} e^{-\frac{4}{5} n} \\
D_{n} = e^{\frac{\pi \, n}{5^{5/4} 3^{1/4}2}} c_{2,n}  \sim \sqrt{\frac{\pi}{2}} \, n^{-1/10} e^{-0.00183 \,n}
\quad \mbox{and} \quad
E_{n} = e^{-\frac{2\pi \, n}{5^{5/4} 3^{1/4}}} c_{2,n}  \sim \sqrt{\frac{\pi}{2}} n^{-1/10} e^{-1.4385 \,n}.
\end{gather*}
\et
\noindent
From (\ref{eq:uniformbdd01a}) it is seen that the $\psi_n$ stay near $\cF^{-1}\chi_{n}$ also piecewise, i.e., inside each interval
given by the successive eigenvalues.

Finally, write $\chi_{n}(y) := C_{n} \Ai_{4}(|y| - \la_{n})$, where $C_{n}:= 1/\Vert \Ai_{4}(|\cdot| - \la_{n}) \Vert_{2}$ for each
$n \in \NN$. Then we have
\bt{thm:diffbound}
For all  $n \in \mathbb N$
$$
\Vert \psi_{n} - {\mathcal F}^{-1}\chi_{n} \Vert_{2} \leqslant \frac{2^{33/20} \, 3^{1/5} \,
5^{1/10}}{\pi^{13/10}} \sqrt{\GA\lk \frac{4}{5} \rk} \, \BT_{n},
$$
where $\BT_{n} := n^{1/5} C_{n} \sqrt{c_{2,n}/c_{1,n}} \asymp
\frac{\pi}{2} \, n^{1/5}\la_{n}^{-1/8} e^{-\frac{2}{5}\la_{n}^{5/4}} \sim \frac{\pi}{2} \, n^{1/10} e^{-\frac{2}{5}n}$.
\et

\bigskip

\section{Proofs}
\subsection{Proof of Proposition \ref{thm:eigenvalexpr01}}\label{sec:lemeigenvalexpr}
We only show the proof for \eqref{eq:eigenvalexpr01b}, by a similar argument \eqref{eq:eigenvalexpr02b} also follows. Let
$\la>0$ and consider the split-up
$
\Ai_{4}(-\la) \widetilde{\Ai}_{4}^{\prime \prime}(-\la) = I(\la) + J(\la),
$
where
$$
I(\la) := \frac{1}{2\pi^{2} i} \lk \int_{-i\infty}^{i\infty} e^{\la t - \frac{t^{5}}{5}} \ud t \rk
\lk \int_{0}^{\infty} u^{2} e^{\la u - \frac{u^{5}}{5}} \ud u \rk
$$
and
$$J(\la) := \frac{1}{\pi^{2}} \lk \int_{0}^{\infty} \cos\lk \frac{t^{5}}{5} - \la t\rk \ud t \rk
\lk \int_{0}^{\infty} u^{2} \sin\lk \frac{u^{5}}{5} - \la u \rk \ud u \rk.
$$
Rewriting $I(\la)$ as a double integral and applying the change of variables $X = t + u$, $Y = t - u$,
we readily obtain that
$
I(\la) = -\frac{1}{16\pi^{2} i} \int_{-i\infty}^{i\infty} e^{\la X + \frac{X^{5}}{20}}
\lk \int_{0}^{\infty}(X+Y)^{2} e^{-\frac{X}{16}(Y^{2} + X^{2})^{2}} \ud Y \rk \ud X.
$
For the second integral a similar straightforward computation gives $J(\lambda) = 0$. Hence
\eq{eq:specfn07}
\Ai_{4}(-\la) \widetilde{\Ai}_{4}^{\prime \prime}(-\la) =
-\frac{1}{16\pi^{2} i} \int_{-i\infty}^{i\infty} e^{\la X + \frac{X^{5}}{20}} \lk \int_{0}^{\infty}(X+Y)^{2}
e^{-\frac{X}{16}(Y^{2} + X^{2})^{2}} \ud Y \rk \ud X.
\en
Next consider similarly
$
\Ai_{4}^{\prime \prime}(-\la) \widetilde{\Ai}_{4}(-\la) = \tilde{I}(\la) + \tilde{J}(\la),
$
where
$$
\tilde{I}(\la) := \frac{1}{2\pi^{2} i} \lk \int_{-i\infty}^{i\infty} t^{2}
e^{\la t - \frac{t^{5}}{5}} \ud t \rk \lk \int_{0}^{\infty} e^{\la u - \frac{u^{5}}{5}} \ud u \rk
$$
and
$$
\tilde{J}(\la) := \frac{1}{\pi^{2}} \lk \int_{0}^{\infty} t^{2} \cos\lk \frac{t^{5}}{5} - \la t\rk \ud t \rk
\lk \int_{0}^{\infty} \sin\lk \frac{u^{5}}{5} - \la u \rk \ud u \rk.
$$
We obtain
$
\tilde{I}(\la) = -\frac{1}{16\pi^{2} i} \int_{-i\infty}^{i\infty} e^{\la X + \frac{X^{5}}{20}}
\lk \int_{0}^{\infty}(X-Y)^{2} e^{-\frac{X}{16}(Y^{2} + X^{2})^{2}} \ud Y \rk \ud X
$
and $\tilde{J}(\la)=0$, which give
\eq{eq:specfn09}
\Ai_{4}^{\prime \prime}(-\la) \widetilde{\Ai}_{4}(-\la)  =
-\frac{1}{16\pi^{2} i} \int_{-i\infty}^{i\infty} e^{\la X + \frac{X^{5}}{20}} \lk \int_{0}^{\infty}(X-Y)^{2}
e^{-\frac{X}{16}(Y^{2} + X^{2})^{2}} \ud Y \rk \ud X.
\en
A combination of \eqref{eq:specfn07}-\eqref{eq:specfn09} gives furthermore
\eq{eq:specfn10}
\Phi_{1}(-\la) = -\frac{1}{4\pi^{2}i} \int_{-i\infty}^{i\infty} X b(X) e^{\la X + \frac{X^{5}}{20}} \ud X
\en
with $b(X) = \int_{0}^{\infty} Y e^{-\frac{X}{16}(Y^{2} + X^{2})^{2}} \ud Y$. Computing the integral, we have
\begin{align*}
\int_{0}^{\infty} Y e^{-\frac{X}{16}(Y^{2} + X^{2})^{2}} \ud Y & =
\frac{1}{2}\int_{0}^{\infty} e^{-\frac{X}{16}(Y + X^{2})^{2}} \ud Y  = X^{-1/2} \GA\lk \frac{1}{2}, \frac{X^{5}}{16} \rk.
\end{align*}
Thus $\Phi_{1}(-\la) = -\frac{1}{4\pi^{2}i} \int_{-i\infty}^{i\infty} X^{1/2} \GA\lk \frac{1}{2}, \frac{X^{5}}{16} \rk
e^{\la X + \frac{X^{5}}{20}} \ud X$, and hence
\begin{align}\label{eq:specfn12}
\Phi_{1}(-\la) &= -\frac{1}{4\pi^{2}} \int_{-\infty}^{\infty} |X|^{1/2} \GA\lk \frac{1}{2}, \frac{iX^{5}}{16} \rk
e^{i\lk \frac{X^{5}}{20} + \la X + \frac{\pi}{4}{\rm sgn}(X)\rk} \ud X \\
&= -\frac{1}{2\pi^{2}} \Re e\int_{0}^{\infty} X^{1/2} \GA\lk \frac{1}{2}, \frac{iX^{5}}{16} \rk e^{i\lk \frac{X^{5}}{20} + \la X
+ \frac{\pi}{4}\rk} \ud X. \non
\end{align}
Finally, we use the identity $\GA\lk \frac{1}{2}, \frac{iX^{5}}{16} \rk = e^{\frac{\pi i}{4}} \lkkk {\rm ci}
\lk \frac{1}{2}, \frac{X^{5}}{16} \rk - {\rm si}\lk \frac{1}{2}, \frac{X^{5}}{16} \rk \rkkk$, see \cite{NIST10}, to complete the
proof for $\Phi_{1}(-\la)$ in \eqref{eq:eigenvalexpr01b}, where
\begin{equation}
\label{fresnel}
{\rm si}(a,z) = \int_{z}^{\infty} t^{a-1} \sin t \ud t \quad \mbox{and} \quad {\rm ci}(a,z) = \int_{z}^{\infty} t^{a-1} \cos t
\ud t \quad \mbox{($\Re e(a)< 1$ and $\Re e(z) > 0$)}
\end{equation}
are the generalised Fresnel sine and cosine integral functions, respectively.

\subsection{Proof of Corollary \ref{cor:asymptrelation}}\label{sec:thmasymptrelation}
Using a steepest descent-argument
we get \eq{eq:fnexp01}
\Ai_{4}(-\la) = \frac{1}{\sqrt{2\pi}}\la^{-3/8}\sin\lk \frac{4}{5}\la^{5/4} + \frac{\pi}{4} \rk \lk 1 + O(\la^{-5/8}) \rk
\en
and
\eq{eq:fnexp02}
\widetilde{\Ai}_{4}(-\la) = \frac{1}{\sqrt{2\pi}}\la^{-3/8} e^{\frac{4}{5}\la^{5/4}} \lk 1 + O(\la^{-5/8}) \rk
\en
for every $\la > \la_{0}$ for some $0 < \la_{0} \leqslant 1/2$. Since $\Ai_{4}(y)$ and $\widetilde{\Ai}_{4}(y)$ are analytic
functions, \eqref{eq:asymprel01} and \eqref{eq:asymprel02} follow by differentiation in \eqref{eq:fnexp01}-\eqref{eq:fnexp02}.
Since these calculations are lengthy but do not use special methods, we leave the details to the reader; alternatively see
\cite{D14}.

\subsection{Proof of Theorem \ref{thm:deviation01}}\label{sec:thmdeviation01}
Before proceeding with the proof of Theorem \ref{thm:deviation01} we will need some auxiliary results. Denote
$\cL := \frac{\ud^{4}}{\ud y^{4}} + |y|$, so that $\lk \cL - \la_{n} \rk \phi_{n}(y) = 0$.
\bl{lem:normestimate01}
For every $n \in \NN$ we have
\eq{eq:normest01}
\Vert \lk \cL - \mu_{n} \rk \Ai_{4}(|\cdot| - \la_{n}) \Vert_{2}^{2} \leqslant \frac{4\sqrt{2}\,\,\GA\lk \frac{4}{5} \rk}
{5^{1/5}\pi}\frac{c_{2,n}}{c_{1,n}}.
\en
\el
\begin{proof}
Since $\mu_{n} \geqslant \la_{n}$ for $n \in \NN$, we have
\begin{align}\label{eq:normestfn01}
| \lk \cL - \mu_{n} \rk \Ai_{4}(|y| - \la_{n}) | &= \frac{1}{c_{1,n}} \big| \lk \cL - \mu_{n} \rk \phi_{n}(y) - c_{2,n}
\lk \cL - \mu_{n} \rk \widetilde{\Ai}_{4}(|y|-\la_{n}) \big| \\
&\leqslant \frac{1}{c_{1,n}} \big| \lk \cL - \la_{n} \rk \phi_{n}(y) - c_{2,n} \lk \cL - \mu_{n} \rk
\widetilde{\Ai}_{4}(|y|-\la_{n}) \big| \non \\
& \leqslant \frac{c_{2,n}}{c_{1,n}} | \lk \cL - \la_{n} \rk \widetilde{\Ai}_{4}(|y|-\la_{n})|. \non
\end{align}
Hence,
\begin{align}\label{eq:normestfn02}
\int_{-\infty}^{\infty} | \lk \cL - \mu_{n} \rk \Ai_{4}(|y|-\la_{n}) |^{2} \ud y
\leqslant  \frac{2 c_{2,n}}{c_{1,n}} \int_{0}^{\infty} \bigg| \lk \frac{\ud^{4}}{\ud y^{4}} + y - \la_{n} \rk
\widetilde{\Ai}_{4}(y-\la_{n})\bigg|^{2} \ud y.
\end{align}
Using that $\widetilde{\Ai}_{4} (z) \asymp \frac{1}{\pi} \int_{0}^{\infty} e^{-zt-\frac{t^{5}}{5}} \ud t$ for every $z \in \RR$,
it follows by \eqref{eq:normestfn02} that
\begin{align}\label{eq:normestfn03}
\Vert &\lk \cL - \mu_{n} \rk \Ai_{4}(|\cdot|-\la_{n}) \Vert_{2}^{2} \leqslant \frac{2 c_{2,n}}{\pi^{2} c_{1,n}} \int_{0}^{\infty} \bigg|
\int_{0}^{\infty} \lk t^{4} + y - \la_{n} \rk e^{-(y-\la_{n})t - \frac{t^{5}}{5}} \ud t \bigg|^{2} \ud y \\
&\leqslant \frac{4 c_{2,n}}{\pi^{2} c_{1,n}} \int_{0}^{\infty} \lk \bigg| \int_{0}^{\infty} t^{4} e^{-(y-\la_{n})t -
\frac{t^{5}}{5}} \ud t \bigg|^{2} + \bigg| \int_{0}^{\infty} (y - \la_{n}) e^{-(y-\la_{n})t - \frac{t^{5}}{5}} \ud t \bigg|^{2} \rk \ud y.
\non
\end{align}
To estimate the inner integral in \eqref{eq:normestfn03} we write $t = (\mu_{n}-y)^{1/4}u$ and $h(u) = \frac{u^{5}}{5} - u$, and apply
a standard calculation to obtain
\eq{eq:normestfn04a}
\Big| \int_{0}^{\infty} t^{4} e^{-(y-\mu_{n})t - \frac{t^{5}}{5}} \ud t \Big| \leqslant
\sqrt{\frac{\pi}{2}} |y-\mu_{n}|^{5/8} e^{-\frac{2\sqrt{2}}{5}|y-\mu_{n}|^{5/4}}.
\en
Similarly, it follows that
\eq{eq:normestfn04b}
\Big| \int_{0}^{\infty} (y - \la_{n}) e^{-(y-\la_{n})t - \frac{t^{5}}{5}} \ud t \Big| \leqslant
\sqrt{\frac{\pi}{2}} |y-\la_{n}|^{5/8} e^{-\frac{2\sqrt{2}}{5}|y-\la_{n}|^{5/4}}.
\en
A combination of \eqref{eq:normestfn04a}-\eqref{eq:normestfn04b} and an application of \eqref{eq:normestfn03} gives that
\begin{align}\label{eq:normestfn05}
\Vert \lk \cL - \mu_{n} \rk \Ai_{4}(|\cdot|-\mu_{n}) \Vert_{2}^{2} &\leqslant \frac{4}{\pi}\frac{c_{2,n}}{c_{1,n}} \int_{0}^{\infty}
|y-\mu_{n}|^{5/4} e^{-\frac{4\sqrt{2}}{5}|y-\mu_{n}|^{5/4}} \ud y \\
&= \frac{4}{\pi}\frac{c_{2,n}}{c_{1,n}} [I_{1}(\mu_{n}) + I_{2}(\mu_{n})]. \non
\end{align}
Here
\begin{align*}
I_{1}(\la_{n}) &:= \int_{0}^{\la_{n}} (\la_{n} - y)^{5/4} e^{-\frac{4\sqrt{2}}{5}(\la_{n}-y)^{5/4}} \ud y =
\frac{1}{\sqrt{2}} \lkkk -\la_{n} e^{-\frac{4\sqrt{2}}{5}\la_{n}^{5/4}} +
\int_{0}^{\la_{n}} e^{-\frac{4\sqrt{2}}{5}(\la_{n}-y)^{5/4}} \ud y\rkkk \\
&\leqslant \frac{1}{\sqrt{2}} \lkkk -\la_{n} e^{-\frac{4\sqrt{2}}{5}\la_{n}^{5/4}} +
\int_{0}^{\infty} e^{-\frac{4\sqrt{2}}{5} z^{5/4}} \ud z\rkkk =
\frac{1}{\sqrt{2}} \lkkk -\la_{n} e^{-\frac{4\sqrt{2}}{5}\la_{n}^{5/4}} + \frac{1}{5^{1/5}}\GA\lk \frac{4}{5} \rk \rkkk,
\end{align*}
and similarly
\begin{align*}
I_{2}(\la_{n}) &:= \int_{\la_{n}}^{\infty} (y - \la_{n})^{5/4} e^{-\frac{4\sqrt{2}}{5}(y -
\la_{n})^{5/4}} \ud y = \frac{1}{\sqrt{2}} \int_{\la_{n}}^{\infty} e^{-\frac{4\sqrt{2}}{5}(y-\la_{n})^{5/4}} \ud y
\leq  %\frac{1}{\sqrt{2}} \int_{0}^{\infty} e^{-\frac{4\sqrt{2}}{5} z^{5/4}} \ud z =
\frac{1}{\sqrt{2}} \frac{1}{5^{1/5}}\GA\lk \frac{4}{5} \rk.
\end{align*}
Thus from \eqref{eq:normestfn05} we obtain
\[
\Vert \lk \cL - \mu_{n} \rk \Ai_{4}(|\cdot|-\la_{n}) \Vert_{2}^{2} \leqslant
\frac{4}{\pi \sqrt{2}}\frac{c_{2,n}}{c_{1,n}} \lkkk \frac{2}{5^{1/5}}\GA\lk \frac{4}{5} \rk -
\la_{n} e^{-\frac{4\sqrt{2}}{5}\la_{n}^{5/4}}  \rkkk.
\]
Since the second term in the bracket is positive, the claim follows.
\end{proof}

Next we derive formulae for the normalization constants explicitly in terms of $n$.
\bl{lem:normconst} For all $n \in \NN$, we have that
\eq{eq:constterms}
c_{1,n} \asymp \sqrt{\frac{\pi}{2}} \la_{n}^{-1/8} \sim \sqrt{\frac{\pi}{2}} \, n^{-1/10} \quad \mbox{and}
\quad c_{2,n} \asymp \sqrt{\frac{\pi}{2}} \la_{n}^{-1/8} e^{-\frac{4}{5}\la_{n}^{5/4}} \sim \sqrt{\frac{\pi}{2}} \,
n^{-1/10} e^{-\frac{4}{5}n}.
\en
\el

\begin{proof}[\bf Proof]
Recall \eqref{eq:normconst01}-\eqref{eq:normconst02} and \eqref{eq:constdem01}-\eqref{eq:constdem03} for $\la > 0$.
From \eqref{eq:fnexp01}-\eqref{eq:fnexp02} we obtain
\begin{gather*}
\Ai_{4}(-\la) = \frac{1}{\sqrt{2\pi}}\la^{-3/8}\sin\lk \frac{4}{5}\la^{5/4} + \frac{\pi}{4} \rk \lk 1 + O(\la^{-5/8}) \rk; \\
\widetilde{\Ai}_{4}(-\la) = \frac{1}{\sqrt{2\pi}}\la^{-3/8} e^{\frac{4}{5}\la^{5/4}} \lk 1 + O(\la^{-5/8}) \rk
\end{gather*}
for every $\la > \la_{0}$ with some $0 < \la_{0} \leqslant 1/2$. Since $\Ai_{4}(y)$ and $\widetilde{\Ai}_{4}(y)$ are analytic
functions, by differentiating and some straightforward calculations we obtain
\begin{gather*}
\LA_{1}(-\la) \asymp \frac{1}{2\pi} \lkkk \sin\lk \frac{4}{5}\la^{5/4} + \frac{\pi}{4} \rk - \cos\lk \frac{4}{5}\la^{5/4}
+ \frac{\pi}{4} \rk \rkkk \la^{-1/2} e^{\frac{4}{5}\la^{5/4}} \\
\LA_{2}(-\la) \asymp \frac{1}{2\pi} \lkkk \sin\lk \frac{4}{5}\la^{5/4} + \frac{\pi}{4} \rk + \cos\lk \frac{4}{5}\la^{5/4}
+ \frac{\pi}{4} \rk \rkkk e^{\frac{4}{5}\la^{5/4}} \\
\LA_{3}(-\la) \asymp \frac{1}{2\pi} \lkkk \cos\lk \frac{4}{5}\la^{5/4} + \frac{\pi}{4} \rk - \sin\lk \frac{4}{5}\la^{5/4}
+ \frac{\pi}{4} \rk \rkkk \la^{1/2} e^{\frac{4}{5}\la^{5/4}}.
\end{gather*}
Using these relations, we have $\la \LA_{1}^{2}(-\la) + \LA_{2}^{2}(-\la) \asymp \frac{1}{2\pi^{2}} e^{\frac{8}{5}\la^{5/4}}$.
In particular,
\[
-\frac{1}{\sqrt{2}} \frac{\widetilde{\Ai}_{4}^{\prime}(-\la_{n})}{\sqrt{\la_{n}\LA_{1}^{2}(-\la_{n}) + \LA_{2}^{2}(-\la_{n})}}
\asymp \sqrt{\frac{\pi}{2}} \la_{n}^{-1/8} \quad \mbox{for every $n = 1, 3, 5, \ldots$}
\]
Similarly, we have
\[
\frac{1}{\sqrt{2}} \frac{\Ai_{4}^{\prime}(-\la_{n})}{\sqrt{\la_{n}\LA_{1}^{2}(-\la_{n}) + \LA_{2}^{2}(-\la_{n})}} \asymp
-\sqrt{\frac{\pi}{2}} \la_{n}^{-1/8} \cos\lk \frac{4}{5}\la_{n}^{5/4} + \frac{\pi}{4} \rk e^{-\frac{4}{5}\la_{n}^{5/4}}
\quad \mbox{for every $n = 1, 3, 5, \ldots$}
\]
We also have that
\begin{align*}
\LA_{2}(-\la) \LA_{3}(-\la) &\asymp \frac{1}{4\pi^{2}} \lkkk \sin^{2}\lk \frac{4}{5}\la^{5/4} + \frac{\pi}{4} \rk -
\cos^{2}\lk \frac{4}{5}\la^{5/4} + \frac{\pi}{4} \rk \rkkk \la^{1/2} e^{\frac{8}{5}\la^{5/4}} \\
&= \frac{1}{4\pi^{2}} \sin\lk \frac{8}{5}\la^{5/4} \rk \la^{1/2} e^{\frac{8}{5}\la^{5/4}}.
\end{align*}
In particular, setting $\la = \lkkk \frac{5}{16}(2n - 1)\pi \rkkk^{4/5}$, we obtain for all $n = 2, 4, 6, \ldots$, that
\[
\LA_{2}(-\la_{n}) \LA_{3}(-\la_{n}) \asymp \frac{1}{4\pi^{2}} \la_{n}^{1/2} e^{\frac{8}{5}\la_{n}^{5/4}}
\]
and thus
\[
\frac{1}{2} \frac{\widetilde{\Ai}_{4}^{\prime \prime}(-\la_{n})}{\sqrt{\LA_{2}(-\la_{n}) \LA_{3}(-\la_{n})}}
\asymp \sqrt{\frac{\pi}{2}} \la_{n}^{-1/8}
\]
and
\[
-\frac{1}{2} \frac{\Ai_{4}^{\prime \prime}(-\la_{n})}{\sqrt{\LA_{2}(-\la_{n}) \LA_{3}(-\la_{n})}} \asymp
\sqrt{\frac{\pi}{2}} \la_{n}^{-1/8} \sin\lk \frac{4}{5}\la_{n}^{5/4} + \frac{\pi}{4} \rk e^{-\frac{4}{5}\la_{n}^{5/4}}.
\]
The fact that $\la_{n} \sim n^{4/5}$ proves the claim.
\end{proof}

Let
\eq{eq:normconstairy01}
C_{n} := 1/\Vert \Ai_{4}(|\cdot| - \la_{n}) \Vert_{2}.
\en
Using the identity $\lk z\Ai_{4}^{2}(z) + 2\Ai_{4}^{\prime}(z) \Ai_{4}^{\prime \prime \prime}(z) -
[\Ai_{4}^{\prime \prime}(z)]^{2} \rk^{\prime} = \Ai_{4}^{2}(z)$, we have that
\begin{align}
\label{eq:normconstairy02}
\Vert \Ai_{4}(|\cdot| - \la_{n}) \Vert_{2} &= \lk \int_{-\infty}^{\infty} \Ai_{4}(|y| - \la_{n}) \ud y \rk^{1/2} =
\lk 2 \int_{-\la_{n}}^{\infty} \Ai_{4}(z) \ud z \rk^{1/2} \non \\
&=\sqrt{2} \lk \la_{n} \Ai_{4}^{2}(-\la_{n}) + [\Ai_{4}^{\prime \prime}(-\la_{n})]^{2} - 2 \Ai_{4}^{\prime}(-\la_{n})
\Ai_{4}^{\prime \prime \prime}(-\la_{n}) \rk^{1/2}.
\end{align}
Recall that $\mu_{n} = -\ap^{\prime}_{n}$ for $n = 1, 3, 5, \ldots$ and $\mu_{n} = -\ap_{n}$ for $n = 2, 4, 6, \ldots$, where
$\ap_{n}$ and $\ap^{\prime}_{n}$ are the negative real zeroes of $\Ai_{4}$ and $\Ai^{\prime}_{4}$, respectively.
Using \eqref{eq:normconstairy01}-\eqref{eq:normconstairy02}, we have that
\eq{eq:gabound01}
C_{n} \asymp \frac{1}{\sqrt{2}} \times \begin{cases} \lk \la_{n}\Ai_{4}^{2}(-\la_{n}) + [\Ai_{4}^{\prime \prime}(-\la_{n})]^{2} \rk^{-1/2}
\quad & n = 2j-1 \\
\lk  [\Ai_{4}^{\prime \prime}(-\la_{n})]^{2} - 2\Ai_{4}^{\prime}(-\la_{n})\Ai_{4}^{\prime \prime \prime}(-\la_{n}) \rk^{-1/2}
\quad & n = 2j\end{cases}
\en
for every $j \in \NN$. We also recall from \eqref{eq:fnexp01} that
\eq{eq:gabound02}
\Ai_{4}(-\la) = \frac{1}{\sqrt{2\pi}} \la^{-3/8}\sin\lk \frac{4}{5}\la^{5/4} + \frac{\pi}{4} \rk \lkkk 1 + O\lk \la^{-5/8} \rk \rkkk
\en
for every $\la > \la_{0}$ given that $0< \la_{0} \leqslant 1/2$. Using that $\Ai_{4}$ is analytic, by differentiation in \eqref{eq:gabound02}
we obtain
\[
\la [\Ai_{4}(-\la)]^{2} + [\Ai_{4}^{\prime \prime}(-\la)]^{2} \asymp \frac{1}{\pi} \la^{1/4}\sin^{2}\lk \frac{4}{5}\la^{5/4} + \frac{\pi}{4} \rk
\]
and
\[
 [\Ai_{4}^{\prime \prime}(-\la)]^{2} - 2\Ai_{4}^{\prime}(-\la)\Ai_{4}^{\prime \prime \prime}(-\la) \asymp \frac{1}{2\pi} \la^{1/4}\lkkk 1+\cos^{2}
 \lk \frac{4}{5}\la^{5/4} + \frac{\pi}{4} \rk\rkkk.
\]
Setting again $\la = [5(2n-1)\pi/16]^{4/5}$ and using \eqref{eq:gabound01}, we obtain
\eq{eq:gabound03}
C_{n} \asymp \sqrt{\frac{\pi}{2}} \la_{n}^{-1/8} \sim \sqrt{\frac{\pi}{2}} \, n^{-1/10}, \quad n \in \NN.
\en
A comparison of \eqref{eq:constterms} and \eqref{eq:gabound03} implies $c_{1,n} \equiv C_{n} := 1/\Vert \Ai_{4}(|\cdot| - \la_{n})\Vert_{2}$.

\begin{proof}[\bf Proof of Theorem \ref{thm:deviation01}]
Since $(\phi_{l})_{l \in \NN}$ forms a complete orthornormal set in $L^{2}(\RR)$, and $\Ai_{4}(|\cdot|-\la_{n})
\in L^{2}(\RR)$, we make the expansion $\Ai_{4}(|y| - \la_{n}) = \sum_{l \geqslant 1} b_{l} \phi_{l}(y)$, $n \in \NN$, where $b_{l}
:= (\Ai_{4}(|\cdot| - \la_{n}), \phi_{l})_{L^2}$, and the Parseval identity
\eq{eq:devtneq02}
\sum_{l \geqslant 1} b_{l}^{2} = \Vert \Ai_{4}(|\cdot|-\la_{n}) \Vert_{2}^{2}
\en
holds. We then have that
\eq{eq:devtneq03}
\Vert (\cL - \mu_{n})\Ai_{4}(|\cdot|-\la_{n}) \Vert_{2}^{2} = \sum_{l \geqslant 1} b_{l}^{2} (\la_{l} - \mu_{n})^{2}.
\en
Let $\la_{q(n)}$ denote an eigenvalue closest to $\mu_{n}$. It follows from \eqref{eq:devtneq03} and an application of \eqref{eq:devtneq02}
that
\eq{eq:devtneq04}
\sum_{l \geqslant 1} b_{l}^{2} (\la_{l} - \mu_{n})^{2} \geqslant (\la_{q(n)} - \mu_{n})^{2}\sum_{l \geqslant 1} b_{l}^{2} = (\la_{q(n)} -
\mu_{n})^{2} \Vert \Ai_{4}(|\cdot| - \la_{n}) \Vert_{2}^{2}.
\en
Using \eqref{eq:devtneq03}-\eqref{eq:devtneq04} and \eqref{eq:normest01}, we obtain
\[
|\la_{q(n)} - \mu_{n}| \leqslant 2\sqrt{\frac{\sqrt{2}\,\,\GA\lk \frac{4}{5} \rk}
{\pi 5^{1/5}}} \frac{\sqrt{c_{2,n}/c_{1,n}}}{\Vert \Ai_{4}(|\cdot| - \la_{n}) \Vert_{2}}.
\]
Since by the above $\Vert \Ai_{4}(|\cdot| - \la_{n}) \Vert^2_{2} = 2\sum_{r=0}^{2}(-1)^{r}\binom{2}{r}\Ai_{4}^{(r)}(-\la_{n})
\Ai_{4}^{(4-r)}(-\la_{n}) \equiv 1/c^2_{1,n}$, we have
\eq{eq:devtneq05}
|\la_{q(n)} - \mu_{n}| \leqslant \frac{2^{5/4}}{5^{1/10}\sqrt{\pi}}\sqrt{\GA\lk \frac{4}{5} \rk} \ga_{n},
\en
with $\ga_{n} := \sqrt{c_{1,n}\, c_{2,n}}$, and using the results of Lemma \ref{lem:normconst}, we have $\ga_{n} \asymp
\sqrt{\frac{\pi}{2}}\,\la_{n}^{-1/8} e^{-\frac{2}{5}\la_{n}^{5/4}} \sim \sqrt{\frac{\pi}{2}} \,n^{-1/10} e^{-\frac{2}{5}n}$.
Take $\del = \frac{5\pi}{4} \lk \frac{\pi}{17} \rk^{5/4} \approx 0.475813$ and let
$n_{0} = \lceil 2 \ln\lk C/\del \rk \rceil +1$,
where $C = \frac{2^{5/4}}{5^{1/10}}\sqrt{\frac{\GA(4/5)}{\pi}} \sup_{n \in \NN} \ga_{n} n^{1/10} e^{2n/5} \approx 0.667134$. We claim that
\begin{equation}
\label{clai}
\la_{q(n)} \in \lk \lk \mu_{n}^{5/4} - \del \rk^{4/5}, \lk \mu_{n}^{5/4} + \del \rk^{4/5} \rk \quad \mbox{for $n \geqslant n_{0}$}.
\end{equation}
By \eqref{eq:devtneq05} we have $|\la_{q(n)} - \mu_{n} | \leqslant C \, n^{-1/10} e^{-2n/5} \leqslant \del$
for $n \geqslant n_{0}$. On the other hand,
\eq{eq:devtneq07}
\frac{(2n-1)\pi}{4} + \max\lkk 0, \frac{(n-3)\pi}{8} \rkk < \mu_{n}^{5/4} - \del \leqslant \mu_{n}^{5/4} + \del \leqslant \frac{5n\pi}{8}
\quad \mbox{for $n \geqslant n_{0}$}.
\en
Using the inequality $\s (\xi - \eta) \xi^{\s - 1} \leqslant \xi^{\s} - \eta^{\s} \leqslant \s (\xi - \eta) \eta^{\s - 1}$ for $0 < \s
\leqslant 1$ and $\xi \geqslant \eta \geqslant 0$, see \cite[(2.15.2)]{bib:halipo52}), and setting
$$
A = \frac{\pi}{C} \lk \frac{\pi}{17} \rk^{5/4} \lk \frac{8}{5\pi} \rk^{1/5} \inf_{n \geqslant n_{0}} n^{-1/10} e^{2n/5} \geqslant 1,
$$
we have that
\begin{align*}
\Big| \lk \mu_{n}^{5/4} \pm \del \rk^{4/5} - \mu_{n} \Big| &> \frac{4 \del}{5} \frac{1}{(\mu_{n}^{5/4} + \del)^{1/5}} \geqslant
\frac{4 \del}{5} \lk \frac{8}{5\pi} \rk^{1/5} n^{-1/5} = \pi \lk \frac{\pi}{17} \rk^{5/4} \lk \frac{8}{5\pi} \rk^{1/5} n^{-1/5} \\
&\geqslant A \, C \, n^{-1/10} e^{-2n/5} \geqslant A |\la_{q(n)} - \mu_{n} | \qquad \mbox{for $n \geqslant n_{0}$}.
\end{align*}
This shows (\ref{clai}). Next by showing that $\mu_{n+1}^{5/4} - \mu_{n}^{5/4} > 2\del$, we conclude that the intervals $\lk (\mu_{n}^{5/4}
- \del)^{4/5}, (\mu_{n}^{5/4} + \del )^{4/5} \rk$ are mutually disjoint, and hence $\la_{q(n)}$ for all $n \geqslant n_{0}$ are
distinct. Using now the inequality $\s (\xi - \eta) \xi^{\s - 1} \geqslant \xi^{\s} - \eta^{\s} \geqslant \s (\xi - \eta) \eta^{\s - 1}$
for $\s \geqslant 1$ and $\xi \geqslant \eta \geqslant 0$, we have
\eq{eq:mutual01}
\mu_{n+1}^{5/4} - \mu_{n}^{5/4} > \frac{5}{4} (\mu_{n+1} - \mu_{n}) \, \mu_{n}^{1/4}.
\en
Take $B = \inf_{n \geqslant n_{0}} (n+1/2)^{-1/5}\mu_{n}^{1/4} \geqslant 1$ and note that (see \eqref{eq:specgaps03} below)
\eq{eq:mutual02}
\mu_{n+1} - \mu_{n} \geqslant \frac{\pi}{2}\lk \frac{8}{15 \pi}\rk^{1/5} \lk n + \frac{1}{2} \rk^{-1/5}, \quad n \in \NN.
\en
A combination of \eqref{eq:mutual02} and \eqref{eq:mutual01} gives
$$
\mu_{n+1}^{5/4} - \mu_{n}^{5/4} > \frac{5\pi}{8}\lk \frac{8}{15 \pi}\rk^{1/5} \lk n + \frac{1}{2} \rk^{-1/5}\mu_{n}^{1/4} \geqslant
\frac{5\pi}{8}\lk \frac{8}{15 \pi}\rk^{1/5} B > 2 \del, \quad \mbox{for all $n \geqslant n_{0}$}.
$$
Furthermore, since $\la_{q(n)}$ is the closest eigenvalue to $\mu_{n}$ so that $\la_{n} \leqslant \la_{q(n)} \leqslant \mu_{n}$ for
$n \geqslant n_{0}$, and given that $\la_{q(n)}^{5/4} < \mu_{n}^{5/4} + \del$ for $n \geqslant n_{0}$, we conclude that $\la_{n} <
\lk \mu_{n}^{5/4} + \del \rk^{4/5}$ for $n \geqslant n_{0}$. Also, we have that
$$
\frac{5n\pi}{8} < \la_{n+1}^{5/4} \leqslant \frac{5(n+1)\pi}{8}, \quad n \in \NN,
$$
which together with \eqref{eq:devtneq07} implies that $\lk \mu_{n}^{5/4} + \del \rk^{4/5} < \la_{n+1}$ for $n \geqslant n_{0}$.
Thus \eq{eq:devtneq08}
\la_{n} < \lk \mu_{n}^{5/4} + \del \rk^{4/5} < \la_{n+1} \quad \mbox{for $n \geqslant n_{0}$}.
\en

To complete, we make use of an idea in \cite[Th. 1]{bib:kwa12}. We claim that there are at most $n_{0} - 1$ eigenvalues
not included in the above class. We use from Section 5.5 below that $\sum_{n \geqslant 1} e^{-\la_{n} t} =
\int_{0}^{\infty} e^{-\la t} \ud N(\la)$, $t > 0$, where $N(\la)$ is the spectral counting function associated with $H$.
Let
$$
L := \lkk l \in \NN \colon l \neq q(n) \,\,\, \mbox{for all $n \geqslant n_{0}$} \rkk.
$$
Note that $L$ is the set of all mismatches $l \neq q(n)$ for $n \geqslant n_{0}$, and hence its complement in $\NN$ is exactly the
set for which $q(n)$ matches $l \in \NN$ for $n \geqslant n_{0}$. Using below $\la_{q(n)} < \lk \mu_{n} + \del \rk^{4/5}$ and the
monotonicity of $e^{-\la t}$, we have
\begin{align*}
\sum_{l \in L} e^{-\la_{l} t}
= \sum_{l \geqslant 1} e^{-\la_{l} t} - \sum_{n \geqslant n_{0}} e^{-\la_{q(n)} t}
%= \int_{0}^{\infty} e^{-\la t} \ud N(\la) - \sum_{n \geqslant n_{0}} e^{-\la_{q(n)} t} \\
\leq \int_{0}^{\infty} e^{-\la t} \ud N(\la) - \sum_{n \geqslant n_{0}} e^{-\lk \mu_{n}^{5/4} + \del \rk^{4/5} t}.
\end{align*}
Since $\int_{\lk \mu_{n_{0}}^{5/4} + \del \rk^{4/5}}^{\infty} e^{-\la t} \ud N(\la) \leqslant
\sum_{n \geqslant n_{0}} e^{-\lk \mu_{n}^{5/4} + \del \rk^{4/5} t}$, we get
\begin{align}\label{eq:devtneq09a}
\sum_{l \in L} e^{-\la_{l} t}  \leq
e^{- \lk \mu_{n_{0}}^{5/4} + \del \rk^{4/5} t} N\lk(\mu_{n_{0}}^{5/4} + \del)^{4/5} \rk + t \int_{0}^{\lk \mu_{n_{0}}^{5/4} +
\del \rk^{4/5}} e^{-\la t} N(\la) \ud \la.
\end{align}
As $N(\la)$ is increasing on $(0, (\mu_{n_{0}}^{5/4} + \del)^{4/5})$, we have
\begin{align}\label{eq:devtneq09b}
\int_{0}^{\lk \mu_{n_{0}}^{5/4} + \del \rk^{4/5}} e^{-\la t} N(\la) \ud \la &\leqslant
 \frac{N\lk (\mu_{n_{0}}^{5/4} + \del)^{4/5} \rk}{t} \lk 1 - e^{-\lk \mu_{n_{0}}^{5/4} + \del \rk^{4/5} t} \rk.
\end{align}
Inserting \eqref{eq:devtneq09b} into \eqref{eq:devtneq09a}, we obtain $\sum_{l \in L} e^{-\la_{l} t} < N\lk (\mu_{n_{0}}^{5/4}
+ \del)^{4/5} \rk$. As $t \searrow 0$, the left-hand side converges to $|L|$, and it follows that $|L| < N\lk (\mu_{n_{0}}^{5/4}
+ \del)^{4/5} \rk$. Using \eqref{eq:devtneq08}, it follows that the right hand side above equals $n$ for $n \geqslant n_{0}$.
Therefore, we have that $|L| < n_{0}$, which proves the claim.

Given that $\la_{n} \leqslant (5n\pi/8)^{4/5}$ for all $n \in \NN$, we have in particular for $l < n_{0}$  that $\la_{l} \leqslant
(5l\pi/8)^{4/5} \leqslant (5(n_{0}-1)\pi/8)^{4/5}$. On the other hand, for $n \geqslant n_{0}$, we have that $\la_{q(n)} > \lk
\mu_{n_{0}}^{5/4} - \del \rk^{4/5} > (5(n-1)\pi/8)^{4/5} \geqslant (5 (n_{0}-1)\pi/8)^{4/5}$. Thus, $L \supseteq \lkk 1, 2, \ldots,
n_{0} - 1 \rkk$. However, $|L| \leqslant n_{0} - 1$, which implies that $L \subseteq \lkk 1, 2, \ldots, n_{0} - 1 \rkk$. Hence they
are actually equal, and therefore we
conclude that $q(n) = n$ for $n \geqslant n_{0}$.

By direct calculation, we have that $n_{0} = 2$ implying that $\la_{1}$ is the only eigenvalue excluded from the above set of
eigenvalues. However, $\la_{1} \leqslant (5\pi/16)^{4/5}$ and since $\mu_{1} \geqslant (\pi/5)^{4/5}$, we obtain
\[
\la_{1} - \mu_{1} \leqslant \lk \frac{5\pi}{16} \rk^{4/5} - \lk \frac{\pi}{5} \rk^{4/5} < \frac{9\pi}{100} \lk \frac{5}{\pi} \rk^{1/5}
\approx 0.310282 < \frac{2^{5/4}}{5^{1/10}}\sqrt{\frac{\GA(4/5)}{\pi}} \, \ga_{1} \approx 0.447194.
\]
Hence we conclude that the result \eqref{eq:devineq01} holds for all $n \geqslant 1$.
\end{proof}

\begin{proof}[\bf Proof of Theorem \ref{thm:specgaps01}]
To derive the upper bound recall from \eqref{eq:eigenvalrel02} that
\eq{eq:specgaps01}
\la_{n+1} - \la_{n} \asymp \lk \frac{5(2n + 1)\pi}{16} \rk^{4/5} - \lk \frac{5(2n - 1)\pi}{16} \rk^{4/5}.
\en
Applying the same inequality to \eqref{eq:specgaps01} as in the previous proof, we obtain
\[
\lk \frac{5(2n + 1)\pi}{16} \rk^{4/5} - \lk \frac{5(2n - 1)\pi}{16} \rk^{4/5} \leqslant \frac{\pi}{2}
\lk \frac{8}{5\pi} \rk^{1/5} \lk n - \frac{1}{2} \rk^{-1/5}.
\]
For the lower bound recall that $\mu_{2j-1} = -\ap^{\prime}_{4,j} \leqslant \lk \frac{5(4j - 3)\pi}{16} \rk^{4/5}$,
$\mu_{2j} = -\ap_{4,j} \geqslant \lk \frac{5(4j - 1)\pi}{16} \rk^{4/5}$, for each $j \in \NN$. The used inequality
gives then again
\eq{eq:specgaps02}
\mu_{2j} - \mu_{2j-1} \geqslant \lk \frac{5(4j - 1)\pi}{16} \rk^{4/5} - \lk \frac{5(4j - 3)\pi}{16} \rk^{4/5}
\geqslant \frac{\pi}{2}\lk \frac{4}{5\pi} \rk^{1/5} \lk j - \frac{1}{4} \rk^{-1/5}.
\en
Setting $n = 2j-1$ in \eqref{eq:specgaps02}, we obtain
\[
\mu_{n+1} - \mu_{n} \geqslant \frac{\pi}{2}\lk \frac{8}{5\pi} \rk^{1/5} \lk n + \frac{1}{2} \rk^{-1/5},
\]
and $n = 2j$  gives
\[
\mu_{n} - \mu_{n-1} \geqslant \frac{\pi}{2}\lk \frac{8}{5\pi} \rk^{1/5} \lk n - \frac{1}{2} \rk^{-1/5} \geqslant
\frac{\pi}{2}\lk \frac{8}{5\pi} \rk^{1/5} \lk n + \frac{1}{2} \rk^{-1/5}.
\]
Thus we conclude that
\eq{eq:specgaps03}
\mu_{n+1} - \mu_{n} \geqslant \frac{\pi}{2}\lk \frac{8}{5\pi} \rk^{1/5} \lk n + \frac{1}{2} \rk^{-1/5} \quad
\mbox{for all $n \in \NN$}.
\en
Since $\la_{n} \asymp \mu_{n}$, using \eqref{eq:specgaps03} we can optimize over $C > 0$ in
\[
\frac{\pi}{2}\lk \frac{8}{5\pi} \rk^{1/5} \lk n + \frac{1}{2} \rk^{-1/5} \geqslant C \lk n - \frac{1}{2} \rk^{-1/5}.
\]
Using $1 \leqslant \frac{n+1/2}{n-1/2} \leqslant 3$ for $n \geq 2$, gives $C = \frac{\pi}{2}\lk \frac{8}{15\pi} \rk^{1/5}$,
and this completes the proof.
\end{proof}

\subsection{Proof of Theorem \ref{thm:traceest}}
Consider the decomposition $\sum_{n=1}^{\infty}e^{-\la_{n} t} = F(t) + G(t)$, where
\[
F(t):= \sum_{j=1}^{\infty}e^{-\la_{2j-1}t} \quad \mbox{and} \quad G(t):= \sum_{j=1}^{\infty}e^{-\la_{2j} t}.
\]
For every $0 \leqslant t \leqslant 1$ it follows from \eqref{eq:eigenvalrel02} that there is a constant $0 < C <1$
such that
\[
e^{-C t }\sum_{j=1}^{\infty}e^{-\lk \frac{5(4j-3)\pi}{16}\rk^{4/5}t } \leqslant F(t) \leqslant
e^{C t } \sum_{j=1}^{\infty}e^{-\lk \frac{5(4j-3)\pi}{16}\rk^{4/5}t }.
\]
Passing from summation to integration gives
\[
e^{-C t }\int_{1}^{\infty} e^{-\lk \frac{5(4x-3)\pi}{16}\rk^{4/5}t} \ud x \leqslant F(t) \leqslant
e^{C t} \int_{1}^{\infty}e^{-\lk \frac{5(4x-3)\pi}{16}\rk^{4/5}t} \ud x.
\]
By the change of variable $y=\lk \frac{5(4x-3)\pi}{16}\rk^{4/5}t$ it furthermore follows that
\[
\int_{1}^{\infty}e^{-\lk \frac{5(4x-3)\pi}{16}\rk^{4/5}t} \ud x =
\frac{1}{\pi t^{5/4}}\int_{\lk\frac{5\pi}{16}\rk^{4/5}t}^{\infty} y^{\frac{1}{4}}e^{-y} \ud y,
\]
which implies that $\lim_{t \downarrow 0}t^{5/4}F(t) = \frac{1}{\pi}\GA\lk \frac{5}{4}\rk$. Similarly, we obtain
$\lim_{t \downarrow 0}t^{5/4} G(t) = \frac{1}{\pi}\GA\lk \frac{5}{4} \rk$. By combining these results
(\ref{eq:traceasympformula01}) follows.

\subsection{Proof of Corollary \ref{thm:weylasymp}.}
Let $m_{l}$ count the multiplicity of $\la_{l}$ for each $l \in \NN$. With $0 < \la_{1} < \la_{2} <  \ldots < \la_{n}
\leqslant R$, for every $R > 0$ and finite, we have
\[
N(\la) =
\begin{cases}
m_{1} + m_{2} + \ldots + m_{l} \quad & \la_{l} \leqslant \la < \la_{l+1}; \,\, l = 1, 2, \ldots, n \\
0 & 0 < \la < \la_{1}.
\end{cases}
\]
Since eigenvalues come with multiplicities, $m_{l} = 1$ for all $l \in \NN$, and we have
\[
N(\la) = \sum_{l=1}^{n} l {\bf 1}_{[\la_{l}, \la_{l+1})}(\la) \quad \mbox{for $0 < \la \leqslant R$},
\]
with $N(\la) = 0$ for $\la \leqslant 0$. Note that $N(\la)$ is right-continuous at each $\la = \la_{n}$, $n \in \NN$,
and jumps by $d_{n} := N(\la_{n}^{+}) - N(\la_{n}^{-}) = 1$, $n \in \NN$. Also, $N(\la)$ is of bounded variation in
$(0,R)$, and the integral $\int_{0}^{R} e^{-\la t} \ud N(\la)$, $t > 0$, exists. Furthermore, taking the limit as $R
\to \infty$, we have that
%\eq{eq:weylproof01}
$\int_{0}^{\infty} e^{-\la t} \ud N(\la) = \lim_{R \to \infty} \int_{0}^{R} e^{-\la t} \ud N(\la)$.
%\en
This limit exists for every $t \in \RR^+$ since by integration-by-parts, $\lim_{R \to \infty} \int_{0}^{R} e^{-\la t}
\ud N(\la) = t \int_{0}^{\infty} e^{-\la t} N(\la) \ud \la < \infty$. Since
\begin{align*}
t \int_{0}^{\infty} e^{-\la t} N(\la) \ud \la &= t \lk \int_{0}^{\la_{1}} +
\sum_{n \geqslant 1} \int_{\la_{n}}^{\la_{n+1}} \rk e^{-\la t} N(\la) \ud \la \\
&= t \sum_{n \geqslant 1} n \int_{\la_{n}}^{\la_{n+1}} e^{-\la t} \ud \la =
\sum_{n \geqslant 1} n \lk e^{-\la_{n} t} - e^{-\la_{n+1} t} \rk.
\end{align*}
By telescopic summation we then obtain
$\int_{0}^{\infty} e^{-\la t} \ud N(\la) = \sum_{n \geqslant 1} e^{-\la_{n} t}$.
By definition, $Z(t)= \sum_{n=1}^{\infty} e^{-\la_{n}t}$, and from Theorem \ref{thm:traceest} we have
\eq{eq:weylproof05}
\int_{0}^{\infty} e^{-\la t} \ud N(\la) = Z(t) \sim \frac{2\GA\lk \frac{5}{4} \rk}{\pi} t^{-5/4} \quad
\mbox{as $t \downarrow 0$}.
\en
A Karamata-Tauberian theorem \cite[Th. 15.3]{bib:kor04} applied to \eqref{eq:weylproof05} gives \eqref{eq:weylasyformula01},
which completes the proof.

\subsection{Proof of Theorem \ref{thm:eigenfns01}}\label{sec:eigenfns01}
We are only required to derive expressions for $c_{1,n}$ and $c_{2,n}$. Take $\phi_{n}(y) = \varphi(y - \la_{n})$ for $y > 0$
and recall from \eqref{eq:diff01} that
\eq{eq:normprf01}
\frac{\ud^{4}}{\ud z^{4}} \varphi(z) = -z \varphi(z).
\en
Consider
$$
\int_{-\infty}^{\infty} \phi_{n}^{2}(y) \ud y = 2\int_{0}^{\infty} \varphi^{2}(y - \la_{n}) \ud y =
2\int_{-\la_{n}}^{\infty} \varphi^{2}(z) \ud z.
$$
Using \eqref{eq:normprf01}, we have the identity
$$
\frac{\ud}{\ud z} \lk z\varphi^{2}(z) + 2 \varphi^{\prime}(z) \varphi^{\prime \prime \prime}(z) -
[\varphi^{\prime \prime}(z)]^{2} \rk = \varphi^{2}(z).
$$
Thus we obtain
\eq{eq:normprf02}
\int_{-\infty}^{\infty} \phi_{n}^{2}(y) \ud y = 2\lk \la_{n}\varphi^{2}(-\la_{n}) +
[\varphi^{\prime \prime}(-\la_{n})]^{2} - 2 \varphi^{\prime}(-\la_{n}) \varphi^{\prime \prime \prime}(-\la_{n})\rk.
\en
From the boundary conditions \eqref{eq:bdrycond01}, we have for each $n = 1, 3, 5, \ldots$ that
\eq{eq:normprf04}
\varphi^{\prime}(-\la_{n}) = \varphi^{\prime \prime \prime}(-\la_{n}) = 0
\en
implying that
\eq{eq:normprf05a}
c_{2} = -c_{1} \frac{\Ai_{4}^{\prime}(-\la_{n})}{\widetilde{\Ai}_{4}^{\prime}(-\la_{n})} \equiv
-c_{1} \frac{\Ai_{4}^{\prime \prime \prime}(-\la_{n})}{\widetilde{\Ai}_{4}^{\prime \prime \prime}(-\la_{n})}
\en
Plancherel's theorem and a combination of \eqref{eq:normprf02}-\eqref{eq:normprf04} yield
\eq{eq:normprf06}
1 = \int_{-\infty}^{\infty} \psi_{n}^{2}(x) \ud x = \int_{-\infty}^{\infty} \phi_{n}^{2}(y) \ud y =
2\lk \la_{n} \varphi^{2}(-\la_{n}) + \varphi^{\prime \prime}(-\la_{n})^{2} \rk \qquad n = 1, 3, 5, \ldots
\en
With \eqref{eq:normprf05a}, we have that
\begin{gather*}
\varphi(-\la_{n}) = c_{1} \Ai_{4}(-\la_{n}) + c_{2} \widetilde{\Ai}_{4}(-\la_{n}) =
-\frac{c_{1}}{\widetilde{\Ai}_{4}^{\prime}(-\la_{n})} \LA_{1}(-\la_{n}) \\
\varphi^{\prime \prime}(-\la_{n}) = c_{1} \Ai_{4}^{\prime \prime}(-\la_{n}) +
c_{2} \widetilde{\Ai}_{4}^{\prime \prime}(-\la_{n}) = \frac{c_{1}}{\widetilde{\Ai}_{4}^{\prime}(-\la_{n})}
\LA_{2}(-\la_{n}).
\end{gather*}
From here we deduce
\begin{gather*}
c_{1,n} \equiv c_{1} =
- \frac{1}{\sqrt{2}} \frac{\widetilde{\Ai}_{4}^{\prime}(-\la_{n})}{\sqrt{\la_{n}\LA_{1}^{2}(-\la_{n}) + \LA_{2}^{2}(-\la_{n})}}  \\
c_{2,n} \equiv c_{2} =
\frac{1}{\sqrt{2}} \frac{\Ai_{4}^{\prime}(-\la_{n})}{\sqrt{\la_{n}\LA_{1}^{2}(-\la_{n}) + \LA_{2}^{2}(-\la_{n})}}
\end{gather*}
for each $n = 1, 3, 5, \ldots$.
To obtain the results for $n = 2, 4, 6, \ldots$, we proceed similarly by using now the boundary conditions \eqref{eq:bdrycond02}
giving
\eq{eq:normprf09}
\varphi(-\la_{n}) = \varphi^{\prime \prime}(-\la_{n}) = 0.
\en
Proceeding as above, we finally obtain
\begin{eqnarray*}
c_{1,n} \equiv c_{1} = \frac{1}{2} \frac{\widetilde{\Ai}_{4}^{\prime \prime}(-\la_{n})}{\sqrt{\LA_{2}(-\la_{n}) \LA_{3}(-\la_{n})}}
\quad \mbox{and} \quad
c_{2,n} \equiv c_{2} = -\frac{1}{2} \frac{\Ai_{4}^{\prime \prime}(-\la_{n})}{\sqrt{\LA_{2}(-\la_{n}) \LA_{3}(-\la_{n})}},
\end{eqnarray*}
for each $n = 2, 4, 6, \ldots$

\subsection{Proof of Theorem \ref{thm:eigenfnasymptotics}}\label{sec:thmeigenfnasymptotics}
We recall \eqref{eq:fouriereigenfns} for each $n = 1, 3, 5, \ldots$, and by taking inverse cosine transform, we have that
\eq{eq:eigenfnasymp02a}
\psi_{2j-1}(x) = \sqrt{\frac{2}{\pi}} \int_{0}^{\infty} \cos(x y) \phi_{2j-1}(y)\ud y \quad \mbox{for each $j \in \NN$},
\en
where
$$
\phi_{2j-1}(y) = c_{1,2j-1} \Ai_{4}(y - \la_{2j-1}) + c_{2,2j-1} \widetilde{\Ai}_{4}(y - \la_{2j-1}), \quad \mbox{for $y > 0$}.
$$
On integrating by parts $(2N+4)$-times in \eqref{eq:eigenfnasymp02a}, together with $\lim_{y \to \infty}
\sin\lk x y + \frac{r\pi}{2}\rk \phi_{2j-1}^{(l)} (y) = 0$ for all $r, l \in \NN_{0}$, we obtain
\begin{align}\label{eq:eigenfnasymp03a}
\int_{0}^{\infty} \cos(x y) \phi_{2j-1}(y) \ud y
&= \sum_{r=0}^{N} (-1)^{r+1} \frac{\phi_{2j-1}^{(2r+1)}(0)}{x^{2r+2}}  + E_{1,N}(x).
\end{align}
Here
\begin{align*}
E_{1,N}(x)
&:= \frac{(-1)^{N+1}}{x^{2N+2}}\int_{0}^{\infty}\cos(x y) \phi_{2j-1}^{(2N+2)}(y) \ud y \\
&= \frac{(-1)^{ N}}{x^{2N+4}}\lk \phi_{2j-1}^{(2N+3)}(0) + \int_{0}^{\infty} \cos(x y) \phi_{2j-1}^{(2N+4)}(y) \ud y \rk,
\end{align*}
for each $N = 2, 3, 4, \ldots$. For sufficiently large $M_{1}, M_{2} > 0$, dependent only on $N$, we have that
$$
|\phi_{2j-1}^{(2N+3)}(0)| \leqslant c_{1,2j-1} |\Ai_{4}^{(2N+3)}(- \la_{2j-1})| + c_{2,2j-1}
|\widetilde{\Ai}_{4}^{(2N+3)}(- \la_{2j-1})| \leqslant M_{1}
$$
and
\begin{align*}
\int_{0}^{\infty} &|\phi_{2j-1}^{(2N+4)}(y)| \ud y \\
&\leqslant c_{1,2j-1} \int_{0}^{\infty}|\Ai_{4}^{(2N+4)}(y - \la_{2j-1})|\ud y + c_{2,2j-1}
\int_{0}^{\infty}|\widetilde{\Ai}_{4}^{(2N+4)}(y - \la_{2j-1})| \ud y \leqslant M_{2}.
\end{align*}
Hence $|E_{1,N}(x)| \leqslant \frac{C_{1,j,N}}{|x|^{2N+4}}$ with some $C_{1,j,N} > 0$. Using the conditions
$\phi_{\la_{2j-1}}^{\prime}(0) = \phi_{\la_{2j-1}}^{\prime \prime \prime}(0) = 0$, $ j = 1, 2, \ldots$, and
differentiating in \eqref{eq:diff01a}, we have $\phi_{2j-1}^{(5)}(y) = -\phi_{2j-1}(y) - (y - \la_{2j-1})
\phi_{2j-1}^{\prime}(y)$, which reduces to $\phi_{2j-1}^{(5)}(0) = -\phi_{2j-1}(0)$ with $\phi_{2j-1}(0) > 0$ for
all $j \in \NN$. Then \eqref{eq:eigenfnasymp01a} follows from \eqref{eq:eigenfnasymp03a} combined with
\eqref{eq:eigenfnasymp02a}. The proof of \eqref{eq:eigenfnasymp01b} follows completely similarly.

\subsection{\bf Proof of Theorem \ref{thm:analytic01}}\label{sec:thmanalytic}
We only show the proof for $n=2j-1$ case, the result is obtained similarly for $n =2j$. For each
$n=2j-1$, $j \in \NN$ we have that
\begin{align}\label{eq:analyticexpprf01a}
\psi_{n}(x) &= \sqrt{\frac{2}{\pi}} \int_{0}^{\infty}\cos(xy) \lkkk c_{1,n}\Ai_{4}(y-\la_{n})
+ c_{2,n}\widetilde{\Ai}_{4}(y-\la_{n})\rkkk \ud y \\
&=
\sqrt{\frac{2}{\pi}} \int_{0}^{\infty} \sum_{r=0}^{\infty}(-1)^{r}\frac{(xy)^{2r}}{(2r)!}\lkkk c_{1,n}\Ai_{4}(y-\la_{n})
+ c_{2,n}\widetilde{\Ai}_{4}(y-\la_{n})\rkkk \ud y. \non
\end{align}
Define
$$
a_{p}(\la_{n}):=\frac{1}{p!} \sqrt{\frac{2}{\pi}} \int_{0}^{\infty}y^{p}\lkkk c_{1,n}\Ai_{4}(y-\la_{n})
+ c_{2,n}\widetilde{\Ai}_{4}(y-\la_{n})\rkkk \ud y, \quad p=2r, \, r \in \NN.
$$
We show that for each $n=2j-1$, $j \in \NN$,
\eq{eq:analyticexpprf01}
\sum_{r=0}^{\infty} a_{2r}(\la_{n}) |x|^{2r} < \infty
\en
for almost all $x \in \RR$; then we can conclude using the monotone and dominated convergence theorems that the
result follows from \eqref{eq:analyticexpprf01a} for each $n=2j-1$, $j \in \NN$, by interchanging summation with
integration. Let
$
R := \lim_{r \rightarrow \infty}\bigg\vert \frac{a_{2r}(\la_{n})}{a_{2r+2}(\la_{n})}\bigg\vert.
$
Since $|\widetilde{\Ai}_{4}(y)| \leqslant |\Ai_{4}(y)|$ for $y > 0$, we have
$$
\Ai_{4}(y) \asymp \frac{1}{2\pi^{1/2}}y^{-\frac{3}{8}}e^{-\frac{4}{5}y^{\frac{5}{4}}}, \quad y>0,
$$
which implies that there exists a constant $C_{n}>0$ such that
$$
| c_{1,n}\Ai_{4}(y-\la_{n}) + c_{2,n}\widetilde{\Ai}_{4}(y-\la_{n}) | < C_{n} e^{-\frac{4}{5}y^{\frac{5}{4}}},
\quad y>0.
$$
Thus we obtain
\begin{align*}
 a_{2r}(\la_{n})
 &\leqslant \frac{C_{n}}{(2r)!}\sqrt{\frac{2}{\pi}} \int_{0}^{\infty}y^{2r}e^{-\frac{4}{5}y^{\frac{5}{4}}} \ud y.
\end{align*}
By the change of variable $z = \frac{4}{5} y^{5/4}$ we have
\eqst
\tilde{a}_{2r} := \frac{1}{(2r)! } \sqrt{\frac{2}{\pi}} \lk \frac{5}{4} \rk^{\frac{8r-1}{5}}
\int_{0}^{\infty} z^{\frac{8r-1}{5}} e^{-z} \ud z = \frac{1}{(2r)!} \sqrt{\frac{2}{\pi}}
\lk \frac{5}{4}\rk^{\frac{8r-1}{5}}\GA\lk\frac{4(2r+1)}{5}\rk .
\enst
Applying Stirling's formula to $\GA(\ap)$, we furthermore have
\begin{align*}
\frac{\tilde{a}_{2r}}{\tilde{a}_{2r+2}} &= (2r+2)(2r+1)\lk \frac{4}{5}\rk^{\frac{8}{5}}
\frac{\GA\lk\frac{4(2r+1)}{5}\rk}{\GA\lk\frac{4(2r+3)}{5}\rk} \\ &\asymp (2r+2)(2r+1)
\lk \frac{4}{5}\rk^{\frac{8}{5}}e^{\frac{8}{5}} \lk \frac{2r+1}{2r+3}\rk^{\frac{4(2r+1)}{5}
- \frac{1}{2}}\lk\frac{5}{4(2r+3)}\rk^{\frac{8}{5}},
\end{align*}
which implies that $R \to\infty$ as $r \to \infty$.

\subsection{\bf Proof of Theorem \ref{thm:diffsupbound}}\label{sec:thmdiffsupbound}
With $\phi_{n}(y) = c_{1,n}\Ai_{4}(|y|-\la_{n}) + c_{2,n}\widetilde{\Ai}_{4}(|y|-\la_{n})$ and $\chi_{n}(y) :=
c_{1,n}\Ai_{4}(|y|-\la_{n})$, we have that
\begin{align}\label{eq:diffsupbd01}
\Big| \psi_{n}(x) - \lk \sF^{-1}\chi_{n} \rk(x) \Big|
&= \frac{1}{\sqrt{2\pi}} \bigg| \int_{-\infty}^{\infty} e^{ixy} \lkkk \phi_{n}(y) - \chi_{n}(y)) \rkkk \ud y \bigg| \non \\
&=\frac{c_{2,n}}{\sqrt{2\pi}} \Big| \int_{-\infty}^{\infty} e^{ixy} \widetilde{\Ai}_{4}(|y|-\la_{n}) \ud y \Big| %\non \\
\leqslant c_{2,n} \sqrt{\frac{2}{\pi}} \int_{0}^{\infty} \big|\widetilde{\Ai}_{4}(y-\la_{n})\big| \ud y.
\end{align}
By a similar procedure leading to \eqref{eq:normestfn04a}, we can also show that
$$
\big|\widetilde{\Ai}_{4}(y-\la_{n})\big| \leqslant
\frac{1}{\sqrt{2\pi}}|y-\la_{n}|^{-3/8}e^{-\frac{2\sqrt{2}}{5}|y-\la_{n}|^{5/4}}.
$$
Using this estimate in \eqref{eq:diffsupbd01}, we have that
\begin{align}\label{eq:diffsupbd02}
\int_{0}^{\infty} \big|\widetilde{\Ai}_{4}(y-\la_{n})\big| \ud y
&\leqslant
\frac{1}{\sqrt{2\pi}}\int_{0}^{\infty} |y-\la_{n}|^{-3/8} e^{-\frac{2\sqrt{2}}{5}|y-\la_{n}|^{5/4}} \ud y \non \\
&=
\frac{1}{\sqrt{2\pi}}\lim_{\EP \to 0} \lkkk \int_{0}^{\la_{n} - \EP} \frac{e^{-\frac{2\sqrt{2}}{5}(\la_{n} - y)^{5/4}}}
{(\la_{n} - y)^{3/8}} \ud y + \int_{\la_{n} + \EP}^{\infty} \frac{e^{-\frac{2\sqrt{2}}{5}(y - \la_{n})^{5/4}}}{(y - \la_{n})^{3/8}}
\ud y \rkkk \non \\
&\leqslant
\frac{2}{\sqrt{2\pi}}\lim_{\EP \to 0} \int_{\EP}^{\infty} z^{-3/8}e^{-\frac{2\sqrt{2}}{5}z^{5/4}} \ud z = \frac{2^{7/4}}{\sqrt{5}}.
\end{align}
The result \eqref{eq:uniformbdd01} then follows by a combination of \eqref{eq:diffsupbd01}-\eqref{eq:diffsupbd02}.

To prove \eqref{eq:uniformbdd01a}, we proceed in a number of steps. Let $p \neq n-1$ and $p \neq n$. We have that
\begin{align}\label{eq:uniformbddeq01}
\int_{\la_{p}}^{\la_{p+1}} \big| \phi_{n}(y) - \chi_{n}(y)\big| \ud y &= c_{2,n} \int_{\la_{p}}^{\la_{p+1}} \big| \widetilde{\Ai}_{4}(y - \la_{n}) \big| \ud y \non \\
&\leqslant \frac{c_{2,n}}{\sqrt{2\pi}} \int_{\la_{p}}^{\la_{p+1}} |y-\la_{n}|^{-3/8}e^{-\frac{2\sqrt{2}}{5}|y-\la_{n}|^{5/4}} \ud y \non \\
&=\frac{c_{2,n}}{\sqrt{2\pi}} \lkkk \int_{0}^{|\la_{p+1} - \la_{n}|} - \int_{0}^{|\la_{p} - \la_{n}|} \rkkk z^{-3/8}e^{-\frac{2\sqrt{2}}{5}z^{5/4}} \ud z =: c_{2,n} I(p,n).
\end{align}
The change of variables $u = \lk \frac{2\sqrt{2}}{5} \rk^{1/2} z^{5/8}$ gives
$$
I(p,n) = \frac{8}{5\sqrt{2\pi}}\lk \frac{5}{2\sqrt{2}} \rk^{1/2} \lkkk \int_{0}^{\lk \frac{2\sqrt{2}}{5} \rk^{1/2}|\la_{p+1} - \la_{n}|^{5/8}} - \int_{0}^{\lk \frac{2\sqrt{2}}{5} \rk^{1/2}|\la_{p} - \la_{n}|^{5/8}} \rkkk e^{-u^{2}} \ud u.
$$
Using the inequality (see \cite[(1)-(3)]{bib:chu55})
\eq{eq:ineqerrorfn}
\frac{\sqrt{\pi}}{2} \lk 1 - e^{-\frac{4}{\pi}\xi^{2}} \rk^{1/2} \leqslant \int_{0}^{\xi} e^{-u^{2}}
\ud u \leqslant \frac{\sqrt{\pi}}{2} \lk 1 - e^{-\frac{4}{\pi}\xi^{2}} \rk^{1/2}, \qquad \xi \geqslant 0,
\en
we get
\begin{align}\label{eq:uniformbddeq02}
I(p,n) &\leqslant \frac{2^{3/4}}{\sqrt{5}} \lkkk \lk 1 - e^{-\frac{8\sqrt{2}}{5\pi}|\la_{p+1} - \la_{n}|^{5/4}} \rk^{1/2}
- \lk 1 - e^{-\frac{2\sqrt{2}}{5}|\la_{p} - \la_{n}|^{5/4}} \rk^{1/2} \rkkk \non \\
&\leqslant \frac{2^{3/4}}{\sqrt{5}} \lk e^{-\frac{2\sqrt{2}}{5}|\la_{p} - \la_{n}|^{5/4}} -
\frac{1}{2}e^{-\frac{8\sqrt{2}}{5\pi}|\la_{p} - \la_{n}|^{5/4}} \rk \leqslant
\frac{2^{3/4}}{\sqrt{5}} e^{-\frac{2\sqrt{2}}{5}|\la_{p} - \la_{n}|^{5/4}}.
\end{align}
Given $n \in \NN$, for $p > n$ we have by \eqref{eq:specgaps01b} that
$$
\la_{n+q} - \la_{n+q-1} \geqslant \frac{\pi}{2} \lk \frac{8}{15\pi} \rk^{1/5} (n+q-1)^{-1/5}, \quad q \in \NN
$$
so that
\begin{align*}
\la_{n+q} - \la_{n} \geqslant
\frac{\pi}{2} \lk \frac{8}{15\pi} \rk^{1/5} q (n+q-1)^{-1/5}.
\end{align*}
For $n \in \NN$, write $p = n+q$ for every $q \in \NN$ giving
$$
\la_{p} - \la_{n} \geqslant\frac{\pi}{2} \lk \frac{8}{15\pi} \rk^{1/5} (p-n) (p-1)^{-1/5}.
$$
Using $(p-n)^{5/4}(p-1)^{-1/4} \geqslant p (1-n/p)^{5/4} \geqslant p(1-\frac{5n}{4p}) = p - \frac{5}{4}n$, it
follows by \eqref{eq:uniformbddeq02} that
\eq{eq:uniformbddeq03}
I(p,n) \leqslant \frac{2^{3/4}}{\sqrt{5}} e^{-\frac{2\sqrt{2}}{5}\lk \frac{\pi}{2} \rk^{5/4}
\lk \frac{8}{15\pi} \rk^{1/4} \lk p - \frac{5}{4}n\rk} \quad \mbox{for each $p \geqslant n+1$}.
\en
Next for $n \in \NN$ and $p < n-1$, again we have by \eqref{eq:specgaps01b} that
$$
\la_{n+1-j} - \la_{n-j} \geqslant \frac{\pi}{2} \lk \frac{8}{15\pi} \rk^{1/5} (n-j)^{-1/5}, \quad 1 \leqslant j \leqslant n-1
$$
giving
\begin{align*}
\la_{n} - \la_{n-q} \geqslant \frac{\pi}{2} \lk \frac{8}{15\pi} \rk^{1/5} q (n-1)^{-1/5}.
\end{align*}
Setting $p = n - q$ for every $2 \leqslant q \leqslant n-1$, we get
$$
\la_{n} - \la_{p} \geqslant \frac{\pi}{2} \lk \frac{8}{15\pi} \rk^{1/5} (n - p) (n - 1)^{-1/5}.
$$
Using $(n - p)^{5/4}(n - 1)^{-1/4} \geqslant n (1-p/n)^{5/4} \geqslant n(1-\frac{5p}{4n}) = n - \frac{5}{4}p$, it follows
from \eqref{eq:uniformbddeq02} that
\eq{eq:uniformbddeq04}
I(p,n) \leqslant \frac{2^{3/4}}{\sqrt{5}} e^{-\frac{2\sqrt{2}}{5}\lk \frac{\pi}{2} \rk^{5/4} \lk \frac{8}{15\pi} \rk^{1/4}
\lk n - \frac{5}{4}p\rk}, \quad 1 \leqslant p \leqslant n-2.
\en
When $p = n$, it directly follows from \eqref{eq:uniformbddeq01} that
\begin{align*}
\int_{\la_{n}}^{\la_{n+1}} \Big| \phi_{n}(y) - \chi_{n}(y)\Big| \ud y &= c_{2,n} \int_{\la_{n}}^{\la_{n+1}}
\big| \widetilde{\Ai}_{4}(y - \la_{n}) \big| \ud y \non \\
&\leqslant
\frac{c_{2,n}}{\sqrt{2\pi}} \int_{\la_{n}}^{\la_{n+1}} |y-\la_{n}|^{-3/8}e^{-\frac{2\sqrt{2}}{5}|y-\la_{n}|^{5/4}} \ud y \non \\
&=
\frac{c_{2,n}}{\sqrt{2\pi}} \int_{0}^{\la_{n+1} - \la_{n}} z^{-3/8}e^{-\frac{2\sqrt{2}}{5}z^{5/4}} \ud z =: c_{2,n} I(n).
\end{align*}
By a similar change of variables $u = \lk \frac{2\sqrt{2}}{5} \rk^{1/2} z^{5/8}$ and a use of \eqref{eq:ineqerrorfn},
we have
\begin{align*}
I(n) &= \frac{2^{3/4}}{\sqrt{5}} \int_{0}^{\lk \frac{2\sqrt{2}}{5} \rk^{1/2}(\la_{n+1} - \la_{n})^{5/8}} e^{-u^{2}} \ud u \non \\
&\leqslant \frac{2^{3/4}}{\sqrt{5}} \lk 1 - e^{-\frac{8\sqrt{2}}{5\pi}(\la_{n+1} - \la_{n})^{5/4}} \rk^{1/2}
\leqslant \frac{2^{3/4}}{\sqrt{5}} \lk 1 - \frac{1}{2}e^{-\frac{8\sqrt{2}}{5\pi}(\la_{n+1} - \la_{n})^{5/4}}\rk.
\end{align*}
From \eqref{eq:specgaps01b} we immediately have that $\la_{n+1} - \la_{n} \geqslant \frac{\pi}{2} \lk \frac{8}{15\pi} \rk^{1/5} n^{-1/5}$,
and hence
\eq{eq:uniformbddeq05}
\int_{\la_{n}}^{\la_{n+1}} \Big| \phi_{n}(y) - \chi_{n}(y)\Big| \ud y \leqslant \frac{2^{3/4}}{\sqrt{5}} \lk 1 -
\frac{1}{2} e^{-\frac{8}{5^{5/4} 3^{1/4}} n^{-1/4}} \rk c_{2,n}, \quad n \in \NN.
\en
Similarly, if $p = n-1$, we have
$$
\int_{\la_{n-1}}^{\la_{n}} \Big| \phi_{n-1}(y) - \chi_{n}(y)\Big| \ud y
\leqslant \frac{c_{2,n}}{\sqrt{2\pi}} \int_{0}^{\la_{n} - \la_{n-1}} z^{-3/8} e^{-\frac{2\sqrt{2}}{5}z^{5/4}} \ud z =: c_{2,n} I(n-1).
$$
We then have that
$$
I(n-1) = \frac{2^{3/4}}{\sqrt{5}} \int_{0}^{\lk \frac{2\sqrt{2}}{5} \rk^{1/2}(\la_{n} - \la_{n-1})^{5/8}} e^{-u^{2}} \ud u
\leqslant \frac{2^{3/4}}{\sqrt{5}} \lk 1 - \frac{1}{2}e^{-\frac{8\sqrt{2}}{5\pi}(\la_{n} - \la_{n-1})^{5/4}}\rk,
$$
and hence
\eq{eq:uniformbddeq06}
\int_{\la_{n-1}}^{\la_{n}} \Big| \phi_{n-1}(y) - \chi_{n}(y)\Big| \ud y \leqslant \frac{2^{3/4}}{\sqrt{5}}
\lk 1 - \frac{1}{2} e^{-\frac{8}{5^{5/4} 3^{1/4}} n^{-1/4}} \rk c_{2,n}, \quad n \in \NN.
\en
Combining \eqref{eq:uniformbddeq01}, \eqref{eq:uniformbddeq03}-\eqref{eq:uniformbddeq06}, and using $1 - e^{-\varrho n^{-1/4}}
\leqslant \varrho n^{-1/4}$ for $0 < \varrho \leqslant 1$ and $n \geqslant 1$, we obtain
\begin{align}\label{eq:uniformbddeq07}
\sum_{p \geqslant 1} &\int_{\la_{p}}^{\la_{p+1}} \big| \phi_{n}(y) - \chi_{n}(y)\big| \ud y \\
&\leqslant \frac{2^{3/4}}{\sqrt{5}} \lk 1 + \frac{8 n^{-1/4}}{5^{5/4} 3^{1/4}} +
\lk \frac{1}{e^{\frac{2\pi}{5^{5/4} 3^{1/4}}} - 1} + \frac{e^{-\frac{\pi}{15^{1/4} 2}}}{e^{\frac{\pi}{15^{1/4} 2}} - 1} \rk
e^{\frac{\pi}{15^{1/4} 2} n} - \frac{e^{\frac{\pi}{15^{1/4} 2}}
e^{\frac{2\pi}{5^{5/4} 3^{1/4}} n}}{\lk e^{\frac{\pi}{15^{1/4} 2}} - 1\rk}\rk c_{2,n}.\non
\end{align}

Finally, consider
\begin{align}\label{eq:uniformbddeq08}
\int_{-\infty}^{\la_{1}} \big| \phi_{n}(y) - \chi_{n}(y)\big| \ud y &= c_{2,n} \int_{-\infty}^{\la_{1}}
\big| \widetilde{\Ai}_{4}(|y| - \la_{n}) \big| \ud y
= c_{2,n} \lk \int_{0}^{\infty} + \int_{0}^{\la_{1}} \rk \big| \widetilde{\Ai}_{4}(y - \la_{n}) \big| \ud y \non \\
&\leqslant
\frac{c_{2,n}}{\sqrt{2\pi}} \lk \int_{0}^{\infty} + \int_{0}^{\la_{1}} \rk
|y-\la_{n}|^{-3/8}e^{-\frac{2\sqrt{2}}{5}|y-\la_{n}|^{5/4}} \ud y \non \\
&=
\frac{c_{2,n}}{\sqrt{2\pi}} \lk \int_{0}^{\infty} + \int_{0}^{\la_{n}} + \int_{\la_{n} - \la_{1}}^{\la_{n}} \rk
z^{-3/8}e^{-\frac{2\sqrt{2}}{5} z^{5/4}} \ud z \non \\
&\leqslant
\frac{2^{3/4}}{\sqrt{5}} \lkkk 1 + 2 \lk 1 - e^{-\frac{2^{7/2}}{5\pi} \la_{n}^{5/4}} \rk^{1/2} -
\lk 1 - e^{-\frac{2^{3/2}}{5} (\la_{n} - \la_{1})^{5/4}} \rk^{1/2} \rkkk c_{2,n} \non \\
&\leqslant
\frac{2^{3/4}}{\sqrt{5}} \lkkk 2 + e^{-\frac{2^{3/2}}{5} (\la_{n} - \la_{1})^{5/4}} - e^{-\frac{2^{7/2}}{5\pi}
\la_{n}^{5/4}} \rkkk c_{2,n} \non \\
&\leqslant
\frac{2^{3/4}}{\sqrt{5}} \lk 2 + e^{-\frac{2^{3/2}}{5} (\la_{n} - \la_{1})^{5/4}} \rk c_{2,n}.
\end{align}
Since $\la_{n} - \la_{1} \geqslant \frac{\pi}{2} \lk \frac{8}{15\pi} \rk^{1/5} n (n-1)^{-1/5} \geqslant
\frac{\pi}{2} \lk \frac{8}{15\pi} \rk^{1/5} n^{4/5}$, it follows from \eqref{eq:uniformbddeq08} that
\eq{eq:uniformbddeq09}
\int_{-\infty}^{\la_{1}} \big| \phi_{n}(y) - \chi_{n}(y)\big| \ud y \leqslant \frac{2^{3/4}}{\sqrt{5}}
\lk 2 + e^{-\frac{2\pi}{5^{5/4} 3^{1/4}} n} \rk c_{2,n}.
\en
Hence \eqref{eq:uniformbdd01a} follows by a combination of \eqref{eq:uniformbddeq07} and \eqref{eq:uniformbddeq09},
which completes the proof.

\subsection{\bf Proof of Theorem \ref{thm:diffbound}}\label{sec:thmdiffbound}
Since $\chi_{n} \in L^{2}(\RR)$, we have
\eq{eq:diffbound01}
\chi_{n}(y) = \sum_{l \geqslant 1} a_{l} \phi_{l}(y)
\en
with $a_{l} := \lka \chi_{n}, \phi_{l} \rka$ and $\sum_{l \geqslant 1} a_{l}^{2} = \Vert \varphi_{n} \Vert_{2}^{2} = 1$. We now consider
\eq{eq:diffbound02}
\Vert (\cL - \la_{n})\chi_{n} \Vert_{2}^{2} = \sum_{l \geqslant 1} a_{l}^{2} (\la_{l} - \la_{n})^{2} = \sum_{l \neq n} a_{l}^{2} (\la_{l} - \la_{n})^{2}.
\en
For $l \neq n$, we have
\[
|\la_{l} - \la_{n}| \geqslant \la_{n+1} - \la_{n} \geqslant \frac{\pi}{2}\lk \frac{8}{15 \pi} \rk^{1/5} n^{-1/5}.
\]
Therefore, from \eqref{eq:diffbound02}, we simply have that
\eq{eq:diffbound03}
\sum_{l \neq n} a_{l}^{2} \leqslant \lk \frac{2}{\pi}\rk^{2} \lk \frac{15\pi}{8} \rk^{2/5} n^{2/5} \Vert (\cL - \la_{n})\chi_{n} \Vert_{2}^{2}.
\en
Furthermore, from \eqref{eq:diffbound01}, we have $\chi_{n}(y) - a_{n} \phi_{n}(y) = \sum_{l \neq n} a_{l} \phi_{l}(y)$ so that using \eqref{eq:diffbound03}, we have
\begin{align}\label{eq:diffbound04}
\Vert \chi_{n} - a_{n} \phi_{n} \Vert_{2}^{2} &= \sum_{l \neq n} a_{l}^{2} \leqslant \frac{4}{\pi^{2}} \lk \frac{15\pi}{8} \rk^{2/5} n^{2/5} \Vert (\cL - \la_{n})\chi_{n} \Vert_{2}^{2} \non \\
&= \frac{4}{\pi^{2}} \lk \frac{15\pi}{8} \rk^{2/5} C_{n}^{2} n^{2/5} \Vert (\cL - \la_{n})\Ai_{4}(\cdot - \la_{n}) \Vert_{2}^{2}.
\end{align}
By similar procedure adopted in the proof of \eqref{eq:normest01}, we can also easily show that
\[
\Vert (\cL - \la_{n})\Ai_{4}(\cdot - \la_{n}) \Vert_{2}^{2} \leqslant \frac{4\sqrt{2}\,\,\GA\lk \frac{4}{5} \rk}{\pi 5^{1/5}} \frac{c_{2,n}}{c_{1,n}}.
\]
Hence, it follows that
\[
\Vert \chi_{n} - a_{n} \phi_{n} \Vert_{2} \leqslant \frac{2}{\pi} \lk \frac{15\pi}{8} \rk^{1/5} \sqrt{\frac{4\sqrt{2}\,\,\GA\lk \frac{4}{5} \rk}{\pi 5^{1/5}}} \, n^{1/5} C_{n} \sqrt{c_{2,n}/c_{1,n}}.
\]
Finally, we observe that $\sum_{n \geqslant 1} a_{n}^{2} = 1$ implying that $a_{n} \leqslant 1$ for each $n \in \NN$. Therefore, we have that $\Vert \chi_{n} - \phi_{n} \Vert_{2} \leqslant \Vert \chi_{n} - a_{n} \phi_{n} \Vert_{2}$, and hence, we complete the proof by recalling Plancherel's theorem.

\bigskip
\noindent
\textbf{Acknowledgments:} The material of this paper is a result of the collaboration of the authors leading to the
PhD thesis \cite{D14} of the first named author. SD gratefully thanks his reviewers and examiners for the feedback,
and Loughborough University for the financial support received during pursuing his project. JL thanks IH\'ES,
Bures-sur-Yvette, for a visiting fellowship, where part of the paper has been written.

\end{document}